\documentclass[12pt]{article}

\usepackage[latin2]{inputenc}
\usepackage{amsmath}
\usepackage{amssymb}
\usepackage{amsfonts}
\usepackage{stmaryrd}
\usepackage{theorem}
\usepackage{graphicx}
\usepackage{epstopdf}
\usepackage[parfill]{parskip}    
\usepackage{hyperref}
\hypersetup{pdftex,colorlinks=true,allcolors=blue}
\usepackage{hypcap}
\usepackage{algorithm}
\usepackage{algpseudocode} 
\usepackage{authblk}
\usepackage{tabto}
\usepackage[titletoc,title]{appendix}

\newtheorem{thm}{Theorem}[section]
\newtheorem{prp}{Proposition}[section]
\newtheorem{lem}{Lemma}[section]
\newtheorem{cor}{Corollary}[section]

\newtheorem{dfn}{Definition}[section]

\newtheorem{prb}{Problem}[section]
\newtheorem{mth}{Method}[section]

\newtheorem{heu}{Heuristic}[section]

\newtheorem{exa}{Example}[section]

\newcommand{\keywords}[1]{\par\noindent{\small{\em Keywords\/}:\ \ #1}}
\newcommand{\mscclass}[1]{\par\noindent{\small{\em MSC class\/}:\ \ #1}}

\newcommand*{\me}{\mathrm{e}}
\newcommand*{\vi}{\mathrm{i}}

\newcommand*{\cf}{\ensuremath{\varphi}}

\newcommand*{\vt}{\ensuremath{\vartheta}}
\newcommand*{\vd}{\ensuremath{\delta}}
\newcommand*{\vl}{\ensuremath{\lambda}}
\newcommand*{\va}{\ensuremath{\alpha}}

\newcommand*{\mes}{\ensuremath{\varnothing}} 
\newcommand*{\N}{\ensuremath{\mathbb{N}}}
\newcommand*{\Z}{\ensuremath{\mathbb{Z}}}
\newcommand*{\Q}{\ensuremath{\mathbb{Q}}}
\newcommand*{\R}{\ensuremath{\mathbb{R}}}
\newcommand*{\C}{\ensuremath{\mathbb{C}}}

\newcommand*{\Conv}{\mathrm{Conv}}

\DeclareMathOperator*{\mamax}{\mathrm{argmax}}
\newcommand*{\argmax}{\mathrm{argmax}}
\newcommand*{\Argmax}{\mathrm{Argmax}}

\newcommand*{\msgn}{\mathrm{sgn}}
\newcommand*{\mint}{\mathrm{int}}
\newcommand*{\mRe}{\mathrm{Re}}
\newcommand*{\mIm}{\mathrm{Im}}

\newcommand*{\Arg}{\mathrm{Arg}}
\newcommand*{\Log}{\mathrm{Log}}

\newcommand*{\co}{\ensuremath{\circ}}

\newcommand{\libeq}{\mathrel{\mathop:}=} 
\newcommand*{\prf}{\textbf{Proof}\ \ }
\newcommand*{\sln}{\textbf{Solution}\ \ }
\newcommand*{\rea}{\textbf{Reasoning}\ \ }
\newcommand*{\sqr}{\ensuremath{\square}}
\newcommand*{\mmod}{\mathrm{mod}}

\newcommand*{\tti}{\ensuremath{\rightarrow\infty}}
\newcommand*{\mydag}{\ensuremath{^{\dag}}}
\newcommand*{\myddag}{\ensuremath{^{\ddag}}}

\newcommand*{\rar}{\ensuremath{\Rightarrow}}

\newcommand*{\lrar}{\ensuremath{\Leftrightarrow}}

\newcommand*{\adr}{\mathrm{adr}}

\newcommand*{\blo}{\mathrm{blo}}
\newcommand*{\ifrm}{\mathrm{ifrm}}
\newcommand*{\Frm}{\mathrm{Frm}}
\newcommand*{\Cyc}{\mathrm{Cyc}}
\newcommand*{\Ext}{\mathrm{Ext}}
\newcommand*{\Inv}{\mathrm{Inv}}
\newcommand*{\Per}{\mathrm{Per}}
\newcommand*{\Eve}{\mathrm{EPer}}
\newcommand*{\Foc}{\mathrm{Foc}}
\newcommand*{\EFoc}{\mathrm{EFoc}}

\newcommand*{\mH}{\mathrm{H}}
\newcommand*{\mM}{\mathrm{M}}
\newcommand*{\T}{\mathrm{T}}

\newcommand*{\mcN}{\ensuremath{\mathcal{N}}}
\newcommand*{\mcT}{\ensuremath{\mathcal{T}}}
\newcommand*{\mcP}{\ensuremath{\mathcal{P}}}
\newcommand*{\mcA}{\ensuremath{\mathcal{A}}}
\newcommand*{\mcAf}{\ensuremath{\mathcal{A}_{fin}}}
\newcommand*{\mcAi}{\ensuremath{\mathcal{A}_\infty}}

\hypersetup{bookmarksnumbered}

\begin{document}

\title{On the Exact Convex Hull of IFS Fractals}
\author{József Vass}
\affil{Faculty of Mathematics\\ University of Waterloo\\ 200 University Avenue West\\ Waterloo, ON, N2L 3G1, Canada\\ jvass@uwaterloo.ca}
\date{February 1, 2018}
\maketitle

\begin{abstract}
\noindent The problem of finding the convex hull of an IFS fractal is relevant in both theoretical and computational settings. Various methods exist that approximate it, but our aim is its exact determination. The finiteness of extremal points is examined a priori from the IFS parameters, revealing some cases when the convex hull problem is solvable. Former results are detailed from the literature, and two new methods are introduced and crystallized for practical applicability -- one more general, the other more efficient. Focal periodicity in the address of extremal points emerges as the central idea.\footnote{Some but not all of these results appeared in the author's doctoral dissertation, and were presented at the 2014 Winter Meeting of the Canadian Mathematical Society, though some updates have been made since then. This research was started in the Fall of 2009, and the results were finalized in February 2015. Some formatting changes were made in 2017. The paper was accepted for publication on Oct. 23, 2017 in the journal Fractals \copyright\ 2018 World Scientific Publishing Company \url{http://www.worldscientific.com/worldscinet/fractals} and it is to appear in \textit{Fractals, Vol. 26, No. 1 (2018) 1850002}.}\\
\ \\
\mscclass{28A80 (primary); 52A10, 52A27 (secondary).}
\keywords{IFS fractals, convex hull, extremal points.}
\end{abstract}

\newpage
\tableofcontents

\newpage
\section{Introduction} \label{s01}

\subsection{The Problem and Former Results} \label{s0102}

The paper investigates the following intentionally vague problem for planar IFS fractals -- defined in Section \ref{s02} -- its vagueness allowing greater freedom of discussion.

\begin{prb} \label{s010104} \textup{(Convex Hull Problem)}\
Devise a general practical method for finding the convex hull of IFS fractals, meaning the exact locations of the extremal points. Determine a priori from the IFS parameters, if the generated fractal will have a finite or infinite number of extremal points.
\end{prb}

Various numerical methods have been devised to find an approximation to the convex hull of IFS fractals, such as \cite{ba00005}, but the goal here is to find the exact convex hull as specified above, which limits the focus and scope of this brief survey.

The following theorem is an intuitive result, which may serve as an inspirational spark for further contemplation of the problem, to one who is already familiar with IFS fractals -- see Section \ref{s02} for an introduction.

\begin{thm} \label{s010101} \textup{(Berger \cite{ic00002})}\
For the extremal points $E_F$ of an IFS fractal $F$ invariant under the Hutchinson operator $\mH(F)=F$, we have that $E_F\subset\mH(E_F)$.
\end{thm}

\begin{cor} \label{s010102} \textup{(Deliu, Geronimo, Shonkwiler \cite{ba00012})}\
For any IFS fractal with a finite number of extremal points, there exists an IFS that generates it, and any extremal address is eventually periodic with respect to this IFS.
\end{cor}

This theorem seems to suggest a recursive property for the set of extremal points, if the above containment is applied repeatedly. If this set is to have a finite cardinality, then the extremal addresses must clearly possess some sort of periodicity, as the corollary states. The authors employ this corollary to resolve the Inverse Problem of IFS Fractals, meaning they devise a theoretical method for finding IFS parameters that generate a given IFS fractal. Their a priori assumption that the fractal has a finite number of extrema seems to be an innocent one. As we will see in Section \ref{s03} it certainly is not, and finiteness will emerge as an intriguing subproblem worthy of our special attention.

The theoretical significance of the Convex Hull Problem \ref{s010104} is underlined by the work of Pearse and Lapidus \cite{ph00005} where the convex hull is an integral tool of their theoretical framework. It is also relevant to applications such as 3D printing \cite{ba00121}. Its relevance for the Fractal-Line Intersection Problem (a.k.a. ``fractal slicing'') is as follows.

\begin{thm} \label{s010103} \textup{\cite{ba00099}}\
An IFS fractal is \textbf{hyperdense} if any hyperplane that intersects its convex hull also intersects the Hutchinson iterate of its convex hull. A hyperplane intersects a hyperdense fractal if and only if it intersects its convex hull. This equivalence holds only if the fractal is hyperdense.
\end{thm}

To highlight a stimulating approach to the Convex Hull Problem, Wang and Strichartz et al. \cite{ba00074, ba00075, ic00004} examine the convex hull of self-affine tiles. They derive the theorem below for the finiteness of extrema for the case of equal IFS factors, which we generalize for unequal factors in Theorems \ref{s030110} and \ref{s030212} via an entirely different approach. Our treatment is restricted to IFS with affine similarities in the complex plane, to allow this generalization. Their approach cannot be generalized to unequal factors, since then the maximization of a linear target function over fractal points in sum form cannot be broken up into individual independent maximizations, as in the equal factor case \cite{ba00074}.

\begin{thm} \label{s010105} \textup{(Strichartz and Wang \cite{ba00074})}\
If the IFS factors are all equal to a certain scaled rotation matrix, then the extremal points of the generated fractal are finite in number if and only if this matrix to some power is a scalar multiple of the identity matrix.
\end{thm}

Mandelbrot and Frame \cite{ba00052} in the process of characterizing ``self-contacting binary trees'' -- with two-map IFS fractal canopies -- determine a periodic address that leads to a particular extremal point. Their result is equivalent to the explicit formula for a periodic point given in Corollary \ref{s020303} of this paper. Their method indicates the relevance of finding the exact convex hull to resolving the Connectedness Problem.

The work most relevant to this paper was done by K{\i}rat and Ko\c{c}yi\u{g}it \cite{ba00109, ma00001, bu00014} focusing essentially on Problem \ref{s010104} for self-affine sets. They devise a constructive terminating algorithm for finding the convex hull of IFS fractals when the factors are equal, and a non-terminating method when they are not. Their methods are presented in a reworked and simplified form in Sections \ref{s0401} and \ref{s0402}, with the important improvement of guaranteed termination.

\subsection{Overview} \label{s0101}

After some preliminary definitions and lemmas in Section \ref{s02}, our investigation begins in Section \ref{s03} with the exclusion of a broad class of IFS fractals -- called ``irrational fractals'' -- for which the cardinality of extremal points is deduced to be infinite, interestingly only dependent on the rotation angles of the IFS maps. The ``rational'' class of fractals on the other hand, contains the broad class of ``fractals of unity'' which as shown in Section \ref{s0302} are guaranteed to have a finite number of extremal points. Furthermore, due to continuity in parameters, fractals of unity can approximate irrational fractals arbitrarily.

Two methods are detailed for finding the convex hull, and contrasted in terms of efficiency and robustness. The first is a ``general method'' that works for any fractal of unity is presented in Section \ref{s0402}, with its shortcomings assessed. The more efficient though less general ``Armadillo Method'' is presented in Section \ref{s0404} for the subclass of ``regular fractals'' which hinges on the linear optimization algorithm of Appendix \ref{s0403}. This method is applied to some examples in Section \ref{s0504}.

The main results are marked with a double dagger \ddag. Ones marked with a dagger \dag\ are either well-known, easy-to-prove, or their proof can be found in the author's PhD thesis \cite{ph00004}.


\subsection{IFS Fractals} \label{s02}

The attractors of iterated function systems -- IFS fractals -- were pioneered by Hutchinson \cite{ba00007} and may be the most elementary fractals possible, occurring in Nature as the Romanesco broccoli. They are the attractors of a finite set of contraction mappings, which when combined and iterated to infinity, converges to an attracting limit set, the IFS fractal itself.

\begin{dfn} \label{s020101}
Let a planar contractive similarity mapping (briefly: \textbf{contraction}) $T:\C\shortrightarrow\C$ be defined for all $z\in\C$ as $T(z)\libeq p+\varphi(z-p)$ where $p\in\C$ is the \textbf{fixed point} of $T$, and $\cf=\vl e^{\vt i}\in\C$ is the \textbf{factor} of $T$, with $\vl\in(0,1)$ the \textbf{contraction factor} of $T$, and $\vt\in(-\pi,\pi]$ the \textbf{rotation angle} of $T$. (Instead of $e^{\vt i}$, rotation matrices can also be used $R\in\R^{d\times d},\ R^{T}\!R=I,\ \det R = 1\ (d\in\N)$.)
\end{dfn}

\begin{dfn} \label{s020102}
Let a planar similarity affine contractive $n$-map \textbf{iterated function system (IFS)} be defined as a finite set of contractions, and denoted as $\mathcal{T}\libeq \{T_1,\ldots,T_n\}$. Denote the \textbf{index set} as $\mcN\libeq \{1,\ldots,n\}$, the respective fixed points as $p_1,\ldots,p_n$ and their set as $\mcP$, and the respective factors as $\cf_1,\ldots,\cf_n$. The fixed points and the factors are called collectively as the \textbf{IFS parameters}. Let the \textbf{Hutchinson operator} be defined as $\mH(S)=\mH_\mcT(S)\libeq \bigcup_{k=1}^n T_k(S)\ (S\subset\C)$.
\end{dfn}

\begin{thm} \label{s020104} \textup{(Hutchinson \cite{ba00007})}
For any IFS $\mcT$ there exists a unique nonempty compact set $F_\mcT\subset\C$ such that $\mH_\mcT(F_\mcT)=F_\mcT$. Furthermore, for any nonempty compact set $S_0\subset\C$, the recursive iteration $S_{n+1}\libeq \mH_\mcT(S_n)$ converges to $F_\mcT$ in the Hausdorff metric.
\end{thm}

\begin{dfn} \label{s020105}
Let the set $F_\mcT$ in the above theorem be called a fractal generated by an IFS $\mcT$ (briefly: \textbf{IFS fractal}; note that multiple IFS can generate the same fractal). Denote $\langle\mathcal{T}\rangle = \langle T_1,\ldots,T_n\rangle \libeq F_\mcT$. If $n=2, 3$ or $n>3$ we say that $F_\mcT$ is a \textbf{bifractal}, \textbf{trifractal}, or \textbf{polyfractal} respectively (see \cite{ba00135} for reasons). If all rotation angles are congruent modulo $2\pi$, then the IFS fractal is said to be \textbf{equiangularly generated} by the IFS $\mathcal{T}$ (briefly: \textbf{equiangular}, if the IFS is understood), and if they are all congruent to zero, then it is a \textbf{Sierpi\'{n}ski fractal}.
\end{dfn}

\begin{figure}[H]
\centering
\includegraphics[width=270pt]{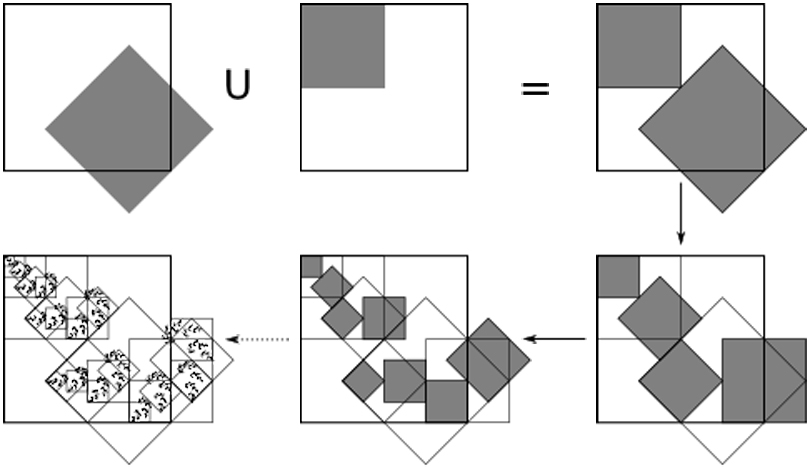}
\caption{Generation of the IFS fractal for $p_1=1+\frac12\vi,\ p_2=\vi,\ \cf_1=\frac{1}{\sqrt{2}}\me^{\frac{\pi}{4}\vi},\ \cf_2=\frac12\me^{0\vi}$ \cite{bw00001}.}
\label{s020106}
\end{figure}

\begin{dfn} \label{s020201}
Let $\mcN^L\libeq \mcN\times\ldots\times\mcN$ be the index set to the $L$-th Cartesian power, and call this $L\in\N$ the \textbf{iteration level}. Then define the \textbf{address set} as
\[ \mcA\libeq \{0\}\cup\bigcup_{L=1}^\infty\mcN^L\cup\mcN^\N. \]
For any $a\in\mcA$ denote its $k$-th \textbf{coordinate} as $a(k),\ k\in\N$. Let its dimension or \textbf{length} be denoted as $|a|\in\N$ so that $a\in\mcN^{|a|}$ and let $|0|\libeq 0$. Define the \textbf{map with address} $a\in\mcA$ acting on any $z\in\C$ as the function composition $T_a(z)\libeq T_{a(1)}\circ\ldots\circ T_{a(|a|)}(z)$ (and the limit if $|a|=\infty$, which exists since the IFS maps are contractive). Let the \textbf{identity map} be $T_0\libeq \mathrm{Id}$. Further denote
\[ \mcAf\libeq \{a\in\mcA: |a|<\infty\},\ \mcAi\libeq \{a\in\mcA: |a|=\infty\}=\mcN^\N. \]
For the weights $w_1,\ldots,w_n\in (0,1)$ let $w_a\libeq w_{a(1)}\cdot\ldots\cdot w_{a(|a|)}$, for the factors $\cf_1,\ldots,\cf_n$ let $\cf_a\libeq \cf_{a(1)}\cdot\ldots\cdot\cf_{a(|a|)}$, and for the angles $\vt_1,\ldots,\vt_n$ let $\vt_a\libeq \vt_{a(1)}+\ldots +\vt_{a(|a|)}$.\\
Let $ab$ denote the \textbf{concatenation} $(a,b)$ for any $a,b\in\mcAf$ so that $T_{ab}=T_a T_b\libeq T_a\circ T_b$ and $a^k\libeq a\ldots a$ $k$-times, $k\in\N$. Let a \textbf{periodic address} be denoted $\bar{a}\libeq aa\ldots$ where $a\in\mcAf$ repeats ad infinitum. Denote the \textbf{inverse} of an address $a\in\mcAf$ with $a^{-1}$ so that $T_{a^{-1}}=T_a^{-1}$ (invertible, since the IFS contraction factors were defined to be non-zero).\\
We say that the address $a\in\mcA$ is a \textbf{truncation} of $b\in\mcA$ if $|a|<\infty$ and there is a $c\in\mcA$ such that $b=ac$, denoted as $a<b$ (note that this includes $a=0$). Furthermore $a\in\mcAf$ is a \textbf{strict truncation} of $b\in\mcA$ if $a<b$ and $a\neq b$ denoted as $a\lneq b$.
\end{dfn}

\begin{prp} \label{s020204}
For any fixed point $p_k\in\mcP$ we have $\langle T_1,\ldots,\T_n\rangle = \lim_{L\rightarrow\infty} \mH^L(\{p_k\}) = \mathrm{Cl}\{T_a(p_k):a\in\mcAf\}\supset \mcP$. This is the \textbf{address generation} of the IFS fractal $F=\langle T_1,\ldots,\T_n\rangle$ from the \textbf{seed} $p_k$. Lastly, $T_a(F),\ |a|=L<\infty$ is said to be a \textbf{subfractal} of $F=\bigcup_{|b|=L} T_b(F)$.\mydag 
\end{prp}

\begin{dfn} \label{s020205}
Let the \textbf{address} $\adr(f)$ of a fractal point $T_a(p_k)$ in the address generation be $a$ if $|a|=\infty$ or $a\bar{k}$ otherwise (if two or more such addresses exist, then take the lexicographically smallest one). We say that a fractal point $f\in F$ is a \textbf{periodic point} if its address $a\in\mcA$ is periodic, meaning it is infinite with a finite repeating part $a=\bar{x},\ x\in\mcAf$. Denote it as $p_x\libeq f$ and note that $T_x(f)=f$. Let the set of all periodic points be denoted as $\Per(F)\libeq \{p_x: x\in\mcA_{fin}\}$. Let the \textbf{cycle} of a finite address $x$ be the set $\Cyc(x)\libeq \{p_{ba}: x=ab\}$.
\end{dfn}

Note that $p_x$ is an abuse of notation that is consistent with the fixed points of the IFS for which $p_k=T_k(p_k),\ k\in\mcN$. Also note that $p_x\in\Cyc(x)$ since $a=0<x$.

\begin{dfn} \label{s020207}
A fractal point $f\in F$ with address $a\in\mcA_\infty$ is \textbf{eventually periodic} (briefly: eventual) if $T_b^{-1}(f)\in\Per(F)$ is periodic for some $b<a$. Let the set of all eventual points be denoted as $\Eve(F)$. (Clearly $\Per(F)\subset\Eve(F)$ with $b=0$.)
\end{dfn}

\begin{lem} \label{s020208}
A fractal point $f\in F$ is eventually periodic iff $\ \exists a,b\in\mcAf:\ f=T_a T_b T_a^{-1}(f)$.\mydag
\end{lem}

\begin{lem} \label{s020302} \textup{(Slope Lemma)}\
$\forall a\in\mcAf,\ z_{1,2}\in\C:\ T_a(z_1)-T_a(z_2) = \cf_a(z_1-z_2)$.\mydag
\end{lem}

\begin{cor} \label{s020303}
A periodic point can be evaluated as $p_x=T_x(0)/(1-\cf_x),\ x\in\mcAf$.\mydag
\end{cor}

\begin{cor} \label{s020304}
The action of a finite map composition is centered at a periodic point as $T_x(z)= p_x+\cf_x(z-p_x)= T_x(0)+\cf_x z$ where $z\in\C,\ x\in\mcAf$.\mydag
\end{cor}

Many beautiful algebraic properties can be shown for periodic points using the Slope Lemma, which we omit for their lack of utility in this paper. Note also that in the first corollary, $0$ is just a part of the identity, and not necessarily a fixed point of the IFS. In closing, two more important lemmas are stated.

\begin{lem} \label{s020301} \textup{(Containment Lemma)}\
If for a nonempty compact set $S\subset\C: \mH(S)\subset S$ then $F\subset S$. Also if $\mH(S)\subset\Conv(S)$ (where $\Conv$ denotes the convex hull), then $F\subset\Conv(S)$. On the other hand, if $F\subset S$ then $F\subset\mH^L(S)$ for any $L\in\N$.\mydag
\end{lem}

\begin{lem} \label{s020305} \textup{(Affine Lemma)}\
For any affine map $M(z)=cz+d\ (c,d,z\in\C)$ we have $M T_k M^{-1}(z) = M(p_k)+\cf_k(z-M(p_k))$ for $k\in\mcN$. Furthermore $M\langle\mcT\rangle = \langle M\mcT M^{-1}\rangle$.\mydag
\end{lem}

\section{On the Cardinality of Vertices} \label{s03}

\subsection{Focal Extremal Points} \label{s0301}

Some preliminary observations are made about the extremal points of polyfractals $F=\langle T_1,\ldots,T_n\rangle$ -- i.e. the vertices of their convex hull -- leading up to the Extremal Theorem at the end of this section. These ideas will prove to be fundamental in the upcoming sections.

\begin{dfn} \label{s030101}
We say that a point $s$ in a compact set $S$ is \textbf{extremal} if it is a vertex of the convex hull $\Conv(S)$, meaning $\nexists s_{1,2}\in S\ (s_1\neq s_2),\ \vl\in(0,1): s=\vl s_1+(1-\vl) s_2$ and call it an \textbf{extremal point}, extremum, or vertex. Denote the set of extremal points as $\Ext(S)$. We say that an \textbf{address is extremal} in an IFS fractal, if the corresponding fractal point is extremal.
\end{dfn}

\begin{lem} \label{s030102} \textup{(Inverse Iteration Lemma)}\
For any extremal point $e\in\Ext(F)$ with address $a=\adr(e)$, taking any truncation $b<a,\ b\in\mcAf$ the \textbf{inverse iterate} $T_b^{-1}(e)$ will also be extremal. Denote the set of all inverse iterates as $\Inv(e)\subset\Ext(F)$.\mydag
\end{lem}

\begin{cor} \label{s030103}
The cycle of a periodic extremum $p_x\in\Ext(F),\ x\in\mcAf$ is also extremal:
\[ \Cyc(x)=\{p_{ba}:\ x=ab\} = \{T_a^{-1}(p_x): a\lneq x\} = \{p_{a^{-1}xa}: a\lneq x\} \subset \Ext(F).\mydag \]
\end{cor}

\begin{cor} \label{s030104}
Inverse iteration of a non-eventual extremal point $e\in\Ext(F)\setminus\Eve(F)$ generates an infinite number of distinct extrema.\mydag
\end{cor}

\begin{dfn} \label{s030105}
Define the \textbf{value} of a finite address $a\in\mcAf\setminus\{0\}$ with respect to an IFS $\mcT$ as the number $\nu(a)=\nu_\mcT(a)\libeq \vt_{a(1)}+\ldots+\vt_{a(|a|)}\ \mmod\ 2\pi\in [0,2\pi)$ and $\nu(0)\libeq 0$, while the value of a set $\mcA_0\subset\mcAf$ as $\nu(\mcA_0)=\{\nu(a): a\in\mcA_0\}$. We say that an address is \textbf{focal} if its value is zero, and \textbf{strictly focal} if it is focal and each index in $\mcN$ occurs as a coordinate. We say that a fractal point is focal, if its address is focal periodic, meaning it is the fixed point of a focal address. Denote the set of all focal points of $F=\langle\mcT\rangle$ as $\Foc(F)=\Foc_\mcT(F)$. We say that a fractal point is \textbf{eventually focal} if it is the finite iterate of a focal point, and denote the set of all eventually focal points as $\EFoc(F)=\EFoc_\mcT(F)$.
\end{dfn}

\begin{thm} \label{s030106}
A periodic extremal point must be focal, meaning $\Ext(F)\cap\Per(F)\subset\Foc(F)$. Furthermore $\mcP\cap\Ext(F)\subset\Foc(F)$, and if $\mcP\subset\Ext(F)$ then $F$ is Sierpi\'{n}ski. Lastly, if $F$ is a Sierpi\'{n}ski fractal then $\Conv(F)=\Conv(\mcP)$ meaning $\Ext(F)\subset\mcP$.
\end{thm}
\noindent
\prf
Suppose that an extremal point $e$ is periodic, meaning there is an $x\in\mcAf$ for which $e=p_x$. Then by Corollary \ref{s020304} we have $T_x(z)=p_x+\cf_x(z-p_x),\ z\in\C$. Suppose indirectly that $\nu(x)\neq 0$ meaning $\cf_x\in\C\setminus\R$. Then the trajectory $t\mapsto T_x^t(z_0)= p_x + \cf_x^t(z_0-p_x)$ traces out a logarithmic spiral for any $z_0\in\C,\ t\in\R$. In the case $\nu(x)\neq\pi$ taking any fractal point $f\in F\setminus\{p_x\}$ the iterates $T_x^l(f)\in F,\ l\in\N\cup\{0\}$ envelop $p_x$ in their convex hull, contradicting that it is extremal. When $\nu(x)=\pi$ the planar segment $[f, T_x(f)]$ contains $p_x$ in its interior, again contradicting its extremality. Thus we must have $\nu(x)=0$, implying by definition that $e=p_x\in\Foc(F)$.\\
The property $\mcP\cap\Ext(F)\subset\Foc(F)$ follows from the previous one, since $\mcP\subset\Per(F)$. Furthermore if $\mcP\subset\Ext(F)$ then $\mcP=\mcP\cap\Ext(F)\subset\Foc(F)\subset\Per(F)\cap\Ext(F)\subset\Foc(F)$ implying that $\nu(k)=0$ meaning $F$ is Sierpi\'{n}ski by definition.\\
Lastly, we show that for a Sierpi\'{n}ski fractal $\Conv(F)=\Conv(\mcP)$ (which is clearly equivalent to $\Ext(F)\subset\mcP$). The $\Conv(\mcP)\subset\Conv(F)$ containment holds since $\mcP\subset F$, and for the reverse containment $\Conv(F)\subset\Conv(\mcP)$ we show that $F\subset\Conv(\mcP)$. All rotation angles are congruent to zero in a Sierpi\'{n}ski fractal, so the IFS contractions are of the form $T_k(z)=(1-\vl_k)p_k + \vl_k z,\ \vl_k\in (0,1)$. Further denoting $F_L\libeq \{T_a(p): |a|=L\}\subset F$ for some $p\in\mcP$ fixed seed, we show that $F_L\subset\Conv(\mcP)$ for all $L\in\N$, which implies that $F\subset\Conv(\mcP)$ as desired. We show this by induction over $L$. For $L=1$ and any $k\in\mcN$ we have $T_k(p)=(1-\vl_k)p_k + \vl_k p\in\Conv(\mcP)$ so $F_1\subset\Conv(\mcP)$. Now suppose that $F_L\subset\Conv(\mcP)$ for some $L\in\N$, and take any $f\in F_L$. Then for any $k\in\mcN$ we have $T_k(f)=(1-\vl_k)p_k + \vl_k f\in\Conv(\mcP)$ thus $F_{L+1}=\bigcup_{k\in\mcN} T_k(F_L)\subset\Conv(\mcP)$. \sqr

\begin{dfn} \label{s040403}
We say that $\tau\in\C\setminus\{0\}$ is a \textbf{target direction} (briefly: target) if our aim is to maximize the linear \textbf{target function} $z\mapsto\langle\tau,z\rangle\libeq (\mRe\ \tau)(\mRe\ z) + (\mIm\ \tau)(\mIm\ z)$ over some compact set $z\in S\subset\C$. The number $\langle\tau,z\rangle$ is called the \textbf{target value} of $z\in\C$. Since IFS fractals are compact, a target function attains its maximum in an extremal point by the Extreme Value Theorem, which we refer to as a \textbf{maximizer} of $\tau$ (or as a maximizing extremal point). Let the set of maximizers be denoted as $\mM(\tau)$ (note that $|\mM(\tau)|\in\{1,2\}$). We call the address of a maximizer, a \textbf{maximizing address} of $\tau$. A truncation of a maximizing address is called a \textbf{maximizing truncation}.
\end{dfn}

\begin{thm} \label{s040404} \textup{(Extremal Theorem; ET)}\
If $a<\adr\ \mM(\tau)$ then $\mM(\tau) = T_a(\mM(\me^{-\nu(a)\vi}\tau))$. Furthermore, if $bx<\adr\ \mM(\tau),\ \nu(x)=0,\ x\neq 0$ then $\mM(\tau)=\{T_b(p_x)\}$. Lastly, if an extremal address begins with $bx\in\mcAf,\ \nu(x)=0,\ x\neq 0$ then that address is $b\bar{x}$.\myddag
\end{thm}
\noindent
\prf
Since $a<\adr\ \mM(\tau)$ the maximum is attained over $T_a(F)$ so we have that
\[ \mM(\tau) = \mamax_{f\in F}\ \langle\tau,f\rangle = \mamax_{f\in T_a(F)}\ \langle\tau,f\rangle = \mamax_{T_a(f_0)\in T_a(F)}\ \langle\tau,\ T_a(f_0)\rangle = \]
\[ = T_a\left(\mamax_{f_0\in F}\ \langle\tau,\ p_a+\cf_a(f_0-p_a)\rangle\right) = T_a\left(\mamax_{f_0\in F}\ \langle\tau,\ \me^{\nu(a)\vi}f_0\rangle\right) = \]
\[ = T_a\left(\mamax_{f_0\in F}\ \langle\me^{-\nu(a)\vi}\tau,\ f_0\rangle\right) = T_a(\mM(\me^{-\nu(a)\vi}\tau)). \]
Thus the first statement holds. Now suppose $bx<\adr\ \mM(\tau),\ \nu(x)=0,\ x\neq 0$. Then
\[ \mM(\tau)\subset T_bT_x(F)\ \rar\ S_0\libeq \mM(\me^{-\nu(b)\vi}\tau)=T_b^{-1}(\mM(\tau))\subset T_x(F)\ \rar\ x<\adr\ \mM(\me^{-\nu(b)\vi}\tau). \]
Letting $\tau'\libeq \me^{-\nu(b)\vi}\tau$ and considering that $|S_0|=|\mM(\tau')|\in\{1,2\}$ we have that $S_0$ is compact, and by the above $S_0=\mM(\tau')=T_x(\mM(\me^{-\nu(x)\vi}\tau'))=T_x(S_0)$ since $\nu(x)=0$. Inductively we see that $S_0=T_x^L(S_0),\ L\in\N$ so taking $L\tti$ we get that $S_0=\{p_x\}$ meaning $\{p_x\}=S_0=\mM(\tau')=\mM(\me^{-\nu(b)\vi}\tau)=T_b^{-1}(\mM(\tau))$ implying $\mM(\tau)=\{T_b(p_x)\}$.\\
The last statement of the theorem follows trivially from the second, considering that a fractal point is extremal iff it maximizes some target direction uniquely. \sqr

\subsection{Rational vs. Irrational Fractals} \label{s0304}

The objective of this section is to show that the cardinality of vertices is infinite for any ``irrational'', while potentially finite for ``rational'' polyfractals $F=\langle T_1,\ldots,T_n\rangle$.

\begin{dfn} \label{s030107}
We say that a polyfractal is \textbf{irrational} with respect to the IFS $\mcT$, if it has no focal points $\Foc_\mcT(F)=\mes$, and \textbf{rational} if it does.
\end{dfn}

\begin{thm} \label{s030110} \textup{(Cardinality of Vertices I)}\
If a polyfractal has a finite number of extrema, then all must be eventual. An irrational polyfractal has no eventual extrema, and has at least a countably infinite number of extremal points.\myddag
\end{thm}
\noindent
\prf
If there were at least one non-eventual extremal point, then by Corollary \ref{s030104} an infinite number of distinct extrema could be generated.\\
Consider an irrational fractal, and suppose indirectly that some extremal point $e\in\Ext(F)$ is eventual, meaning $\exists b,x\in\mcAf:\ e=T_b(p_x)$. Then by Lemma \ref{s030102} the inverse iterate $p_x=T_b^{-1}(e)\in\Ext(F)$ so $p_x\in\Ext(F)\cap\Per(F)\neq\mes$. By the irrationality of $F$ however $\Foc(F)=\mes$ so by Theorem \ref{s030106} we get the empty intersection $\Ext(F)\cap\Per(F)=\mes$ which is a contradiction.\\
Taking any extremal point, it must be non-eventual by the above, so Corollary \ref{s030104} implies that inverse iterating this extremal point will generate a countably infinite number of distinct extrema. \sqr

The question remains whether the converse is true, that is, do rational polyfractals have a finite number of extrema? We show this for the broad subclass of ``fractals of unity'' in Section \ref{s0302}. But first, we reason why it is sufficient to consider only rational fractals in practice.

Due to the lack of extremal focality in an irrational fractal, the extrema tend to infinitesimally cluster around the vertices of the convex hull (by Theorem \ref{s030110}). Thus the corners are ``rounded off self-similarly'' in an infinite number of extrema, so finding the exact convex hull of an irrational fractal seems hopelessly difficult. We may however turn our efforts to finding the potentially finite convex hull of rational fractals, and fortunately by the theorem below, they approximate the irrational case via the subclass of ``fractals of unity''.

\begin{thm} \label{s030301} \textup{(Continuity in Angular Parameters)}\
Keeping the rotation angles variable as a vector $\vt\libeq (\vt_1,\ldots,\vt_n)$ while the other parameters of the IFS constant, and denoting the resulting attractor as $F[\vt]$, the map $\vt\mapsto F[\vt]$ is continuous in the following domain and range metrics respectively
\[ \vd_\infty(\vt,\vt')\libeq \max_{k\in\mcN} |\vt_k-\vt_k'|\ \ and\ \ d_\infty(F,F')\libeq \sup_{a\in\mcA} |T_a(p)-T_a'(p)| \]
where $p\in\mcP$ is any seed.\mydag
\end{thm}

\subsection{Fractals of Unity} \label{s0302}

Since only rational fractals can have a finite number of extrema -- and thus a determinable convex hull -- we focus on the rational subclass of ``fractals of unity'', aiming to show that they indeed have a finite number of extrema at the end of this section. Theorem \ref{s030301} implies that ir/rational fractals can be approximated continuously by fractals of unity.

\begin{prp} \label{s030201}
If all normalized factors $\frac{\cf_k}{|\cf_k|}=\me^{\vt_k\vi},\ k\in\mcN$ of the IFS are roots of unity, meaning $\frac{\vt_k}{2\pi}\in\Q$, then the fractal is rational. We refer to such a case as a \textbf{fractal of unity}.\mydag
\end{prp}

\begin{lem} \label{s030202}
If the rotation angles of the IFS are of the form $\vt_k=2\pi N_k/M$ where $M\in\N,\ N_k\in\{0,\ldots,M-1\},\ k\in\mcN=\{1,\ldots,n\}$ then the number of possible address values is $|\nu(\mcAf)| = M/\gcd(N_1,\ldots,N_n,M)$.
\end{lem}
\noindent
\prf
\[ \nu(\mcAf) = \left\{\frac{2\pi}{M}(k_1 N_1+\ldots +k_n N_n)\ \mmod\ 2\pi:\ k_j\in \{0\}\cup\N,\ j\in\mcN\right\} = \]
\[ = \frac{2\pi}{M}\left\{(k_1 N_1+\ldots +k_n N_n)\ \mmod\ M:\ k_j\in \{0\}\cup\N,\ j\in\mcN\right\}. \]
B\'{e}zout's Lemma implies that the possible values of $(k_1 N_1+\ldots +k_n N_n)\ \mmod\ M$ are limited to $DN\ \mmod\ M\ (N\in\N)$ where $D=\gcd(N_1,\ldots,N_n)$. Since $|\{DN\ \mmod\ M:\ N\in\N\}| = M/\gcd(D,M)$ we have that $|\nu(\mcAf)|=M/\gcd(N_1,\ldots,N_n,M)$. \sqr

\begin{lem} \label{s030203}
In a fractal of unity, any infinite address $a\in\mcAi$ begins with a truncation $bx<a$ for which $b,x\in\mcAf,\ \nu(x)=0,\ x\neq 0$.
\end{lem}
\noindent
\prf
By the above Lemma \ref{s030202}, the truncations of an infinite address $a\in\mcAi$ cannot all be distinct in value, since that would require an infinite number of possible values. Thus
\[ \exists b,x\in\mcAf:\ bx<a,\ x\neq 0,\ \nu(b)=\nu(bx)=\nu(b)+\nu(x) \]
implying that $\nu(x)=0$. \sqr

\begin{thm} \label{s030204}
In a fractal of unity, all extremal points are eventually focal.
\end{thm}
\noindent
\prf
Take any $e\in\Ext(F)$. Then by the above Lemma \ref{s030203} we have $\exists b,x\in\mcAf\ (\nu(x)=0,\ x\neq 0):\ bx<\adr(e)$. The Extremal Theorem \ref{s040404} gives $e=T_b(p_x)\in\EFoc(F)$. \sqr

\begin{cor} \label{s030205}
A fractal of unity has at least one focal extremal point.\mydag
\end{cor}

A finite number of extrema would imply that all are eventual by Theorem \ref{s030110}, thus Theorem \ref{s030204} is a weaker statement that hints at the possibility of finiteness of extrema, which will soon be shown. The existence of a focal periodic extremal point by the corollary, will be exploited in Sections \ref{s0404} and \ref{s05} to generate the entire set of extremal points. We proceed to introducing the concept of ``irreducibility'' of extremal addresses, which by the above theorem are necessarily eventually focal for fractals of unity. This concept will prove to be critical for showing the finiteness of extrema.

\begin{dfn} \label{s030206}
Define the \textbf{forms} of an eventually focal fractal point $f\in F$ with respect to the IFS $\mcT\ (F=\langle\mcT\rangle)$, as the set
\[ \Frm(f)=\Frm_\mcT(f)\libeq\left\{(b,x)\in\mcAf^2:\ \adr(f)=b\bar{x},\ \nu_\mcT(x)=0,\ x\neq 0\right\} \]
and note that $(b,x)$ also conveniently represents the concatenation $bx\in\mcAf$. A form $(b,x)\in\Frm(f)$ is \textbf{irreducible} if there is no shorter form that represents $f$, meaning $\nexists (c,y)\in\Frm(f): |cy|<|bx|$ ($\lrar\ cy\lneq bx$). If a shorter $(c,y)$ form does exist, it is called a \textbf{reduction} of $(b,x)$. An eventually focal point $f=T_b(p_x)\in F$ is said to be \textbf{in irreducible form} if $(b,x)\in\Frm(f)$ is irreducible (note that $b=0$ is permitted).
\end{dfn}

This definition becomes relevant in light of the Extremal Theorem \ref{s040404} (ET), which is essentially an ``inside-out blow-up property'' for extremal points $e\in\Ext(F)$. Regardless of what comes after the reduction $cy < \adr(e)$, the ET implies that $\adr(e)=c\bar{y}$ meaning $e=T_c(p_y)$. Irreducibility implies that no such reduction or simplification of the $\adr(e)=b\bar{x}$ representation is possible.

\begin{thm} \label{s030207}
Any extremal point $e\in\Ext(F)$ in a fractal of unity has a unique irreducible form, and it is a reduction of any other form representing $e$. Denote it as $\ifrm(e)$.
\end{thm}
\noindent
\prf
Clearly $\forall e\in\Ext(F): \Frm(e)\neq\mes$ by Theorem \ref{s030204}. We first prove the existence of the irreducible form by infinite descent. Suppose indirectly that $e\in\Ext(F),\ a\libeq \adr(e)$ has no irreducible form, meaning any form $(b,x)\in\Frm(e)$ has a reduction $(c,y)\in\Frm(e)$. Clearly $\Frm(e)\neq\mes$ by Theorem \ref{s030204}, so starting with a form it cannot be reduced ad infinitum, since any form must have at least a length of one. So for the address $a=\adr(e)\in\mcAi$ the reductions lead inductively to $(0,a(1))\in\Frm(e)$ which is irreducible, leading to a contradiction of the assumption that $\nexists\ifrm(e)$.\\
We proceed to showing the uniqueness of an irreducible form. Suppose indirectly that there are two distinct irreducible forms $(b,x), (c,y)\in\Frm(e)$. Then $|bx|=|cy|$ otherwise one would be a reduction of the other. Since both $bx, cy < \adr(e)$ we have that $bx=cy$. If $b=c$ then $x=y$ leading to a contradiction, so without the restriction of generality, we may suppose that $b\lneq c$ meaning $\exists z\neq 0:\ c=bz$. Then $\nu(b)=\nu(bx)=\nu(cy)=\nu(c)=\nu(b)+\nu(z)$ implying $\nu(z)=0$, therefore $bz=c\lneq cy < a$ so by the ET $(b,z)\in\Frm(e)$ and it is a reduction of $(c,y)$, contradicting the assumption that the latter is irreducible.\\
Lastly, the irreducible form is shorter than any other form, since if there were another form of the same length, then that too would need to be irreducible (since no shorter form exists), which would contradict the uniqueness of the irreducible form. So it must be a truncation and thus a reduction of any other form. \sqr

\begin{lem} \label{s030208}
If $e=T_b(p_x)\in\Ext(F)$ is in irreducible form, then so is $p_x\in\Inv(e)\subset\Ext(F)$.
\end{lem}
\noindent
\prf
Suppose indirectly that $\exists (c,y)\in\Frm(p_x): |cy|<|0x|=|x|$. Then $e=T_b(p_x)=T_{bc}(p_y)$ and $|bcy|=|b|+|cy|<|b|+|x|=|bx|$ implying that $(bc,y)\in\Frm(e)$ and that it is a reduction of $(b,x)$, contradicting that the latter form is irreducible. \sqr

\begin{lem} \label{s030209}
If $e=T_b(p_x)\in\Ext(F)$ is in irreducible form, then $|\Inv(e)|=|b|+|x|$ where $\Inv(e)=\{T_c^{-1}(e): c\lneq b\}\cup\Cyc(x)$ and this union is disjoint.
\end{lem}
\noindent
\prf
First we show that if $p_x\in\Ext(F)\cap\Foc(F)$ is in irreducible form, then $|\Cyc(x)|=|x|$. Clearly $|\Cyc(x)|\leq |x|$ and suppose indirectly that $|\Cyc(x)|<|x|$. Then the cycle must have two identical elements $\exists b,c\ (b\lneq c\lneq x):\ T_b^{-1}(p_x)=T_c^{-1}(p_x)$. Then $\exists y\neq 0:\ c=by$ so $T_c^{-1}=T_y^{-1}T_b^{-1}$ implying that $T_y(T_c^{-1}(p_x)) = T_b^{-1}(p_x) = T_c^{-1}(p_x)$ meaning $p_y= T_c^{-1}(p_x)\in\Cyc(x)\subset\Ext(F)$ which by Theorem \ref{s030106} implies $\nu(y)=0$. Thus we have $by=c\lneq x,\ \nu(y)=0,\ y\neq 0$ so $|by|<|x|$ which by the ET implies that $p_x=T_b(p_y)$ and so $(b,y)\in\Frm(p_x)$ is a reduction of $(0,x)$, contradicting that the latter is irreducible. Thus $|\Cyc(x)|=|x|$.\\
For the general case, take an $e=T_b(p_x)\in\Ext(F)$ in irreducible form, and first note that $\Inv(e)=\{T_c^{-1}(e): c\lneq b\}\cup\Cyc(x)$. Since by Lemma \ref{s030208} the extremal point $p_x$ is also in irreducible form, we have by the above that $|\Cyc(x)|=|x|$, so we only need to show that $|\{T_c^{-1}(e): c\lneq b\}|=|b|$ and $\{T_c^{-1}(e): c\lneq b\}\cap\Cyc(x)=\mes$ to arrive at $|\Inv(e)| = |b|+|x|$.\\
Clearly $|\{T_c^{-1}(e): c\lneq b\}|\leq |b|$ holds, and suppose indirectly that $|\{T_c^{-1}(e): c\lneq b\}|<|b|$. Then $\exists c,d\lneq b\ (c\neq d):\ T_c^{-1}(e)=T_d^{-1}(e)$ and we can assume that $c\lneq d$ without restricting generality. Thus $\exists y\neq 0: d=cy$ so $T_y(T_d^{-1}(e))=T_c^{-1}(e)=T_d^{-1}(e)$ implying that $p_y=T_d^{-1}(e)$ but since $e\in\Ext(F)$ and $d\lneq b<\adr(e)$ we have that $p_y\in\Inv(e)\subset\Ext(F)\cap\Per(F)\subset\Foc(F)$ by Theorem \ref{s030106}, so we can conclude that $\nu(y)=0$. Since $p_y=T_d^{-1}(e)=T_c^{-1}(e)$ we have that $T_c(p_y)=e,\ \nu(y)=0$ implying $(c,y)\in\Frm(e)$ and $|cy|=|d|<|b|<|bx|$, meaning that $(c,y)$ is a reduction of $(b,x)$ contradicting that the latter is irreducible. So we have that $|\{T_c^{-1}(e): c\lneq b\}|=|b|$.\\
Lastly, we show that $\{T_c^{-1}(e): c\lneq b\}\cap\Cyc(x)=\mes$. Suppose indirectly that $\exists c\lneq b, d\lneq x:\ e_0\libeq T_c^{-1}(T_b(p_x))= T_d^{-1}(p_x)$. Then clearly $\exists y,z\neq 0:\ b=cy, x=dz$ so we have $T_b^{-1} = T_y^{-1}T_c^{-1}$ which implies $T_d^{-1}(p_x)=T_c^{-1}(T_b(p_x))= T_y(p_x)$ and thus $T_{dy}(p_x)=p_x$. Therefore $p_{dy}=p_x\in\Ext(F)$ which by Theorem \ref{s030106} implies $\nu(yd)=\nu(dy)=0$. Furthermore $yd\neq 0$ since $y\neq 0$, and also $cyd=bd\lneq bx$ so by the ET we have $e=T_b(p_x)=T_c(p_{yd})$ meaning $(c,yd)\in\Frm(e)$ is a reduction of $(b,x)$ which contradicts that the latter is irreducible. \sqr

\begin{lem} \label{s030210}
If an eventually focal extremal point $e=T_b(p_x)\in\Ext(F)$ is in irreducible form, then so are its inverse iterates $\Inv(e)=\{T_c^{-1}(e): c\lneq b\}\cup\Cyc(x)$.
\end{lem}
\noindent
\prf
First we show this for focal extrema, meaning if $e=p_x\in\Ext(F)\cap\Foc(F)$ is in irreducible form, then for any $a\lneq x$ the form $(0,a^{-1}xa)\in\Frm(p_{a^{-1}xa})$ is irreducible. Suppose indirectly that for some $a\lneq x$ it is reducible, meaning $\exists (b,y)\in\Frm(p_{a^{-1}xa}):\ |by|<|a^{-1}xa|=|x|$ and we may assume that this $(b,y)$ is the irreducible form of $p_{a^{-1}xa}$ by Theorem \ref{s030207}. Then by the above Lemma \ref{s030209} we have that $|by|=|\Inv(T_b(p_y))|=|\Inv(p_{a^{-1}xa})|=|\Cyc(a^{-1}xa)|=|\Cyc(x)|=|x|$ where the latter equality holds since $p_x$ was assumed to be in irreducible form. Therefore $|by|=|x|$ which contradicts $|by|<|x|$ above.\\
For the general case, take an $e=T_b(p_x)\in\Ext(F)$ in irreducible form. Lemma \ref{s030208} implies that $p_x$ is also in irreducible form, so by the above the elements of $\Cyc(x)$ are as well. Thus we only need to show that $(c^{-1}b,x)\in\Frm(T_c^{-1}(e))$ is irreducible for any $c\lneq b$. Clearly $\exists d\neq 0: b=cd$ meaning $d=c^{-1}b$, so we need to show that $(d,x)\in\Frm(T_c^{-1}(e))$ is irreducible. Suppose indirectly that a reduction exists $\exists (a,y)\in\Frm(T_c^{-1}(e)):\ T_a(p_y)=T_d(p_x),\ |ay|<|dx|$. Then $e=T_b(p_x)=T_cT_d(p_x)=T_{ca}(p_y)$ and $|cay|=|c|+|ay|<|c|+|dx|=|cdx|=|cd|+|x|=|b|+|x|=|bx|$ implying that $(ca,y)\in\Frm(e)$ is a reduction of $(b,x)\in\Frm(e)$ contradicting that the latter form is irreducible. \sqr

These theorems show that Definition \ref{s030206} is not merely an intuitive definition, but the proper way to define irreducibility, in order to arrive at the expected properties above.

\begin{thm} \label{s030211}
A focal extremal point cannot have an irreducible form longer than $|\nu(\mcAf)|$, and an eventually focal extremal point cannot have an irreducible form longer than $2|\nu(\mcAf)|$.
\end{thm}
\noindent
\prf
First consider indirectly a focal $e=p_x\in\Ext(F)$ in irreducible form, for which $|x|>|\nu(\mcAf)|$. Then the truncations of $x$ cannot all be different in value, so $\exists b,c\lneq x\ (b\lneq c): \nu(b)=\nu(c)$. Thus $\exists y\neq 0: c=by$ and $\nu(b)=\nu(c)=\nu(b)+\nu(y)$ which implies $\nu(y)=0$, so since $by=c<x<\adr(e)$ by the ET $e=T_b(p_y)$ meaning $(b,y)\in\Frm(e)$ and $|by|=|c|<|x|$ which contradicts that $(0,x)\in\Frm(e)$ is irreducible.\\
Now consider indirectly an eventually focal $e=T_b(p_x)\in\Ext(F)$ in irreducible form, for which $|bx|>2|\nu(\mcAf)|$. Then either $|x|>|\nu(\mcAf)|$ or $|b|>|\nu(\mcAf)|$ must hold (if neither held, then $|bx|\leq 2|\nu(\mcAf)|$). But by Lemma \ref{s030208} we know that $p_x$ is also in irreducible form, so by the above $|x|\leq |\nu(\mcAf)|$ thus necessarily $|b|>|\nu(\mcAf)|$ must hold. Then the truncations of $b$ cannot all be different in value, so $\exists c,d\lneq b\ (c\lneq d): \nu(c)=\nu(d)$. Thus $\exists y\neq 0: d=cy$ and $\nu(c)=\nu(d)=\nu(c)+\nu(y)$ which implies $\nu(y)=0$, so since $cy=d<b<\adr(e)$ by the ET $e=T_c(p_y)$ meaning $(c,y)\in\Frm(e)$ and $|cy|=|d|<|b|<|bx|$ contradicting that $(b,x)\in\Frm(e)$ is irreducible. \sqr

\begin{thm} \label{s030212} \textup{(Cardinality of Vertices II)}\
A fractal of unity $F=\langle T_1,\ldots,T_n\rangle$ has a finite number of extremal points, specifically $|\Ext(F)|\leq n^{2|\nu(\mcAf)|}$.\myddag
\end{thm}
\noindent
\prf
By Theorem \ref{s030204} all extremal points are eventually focal, and by Theorem \ref{s030207} all can be written in irreducible form. By Theorem \ref{s030211} the irreducible form of an eventually focal extremal point cannot be longer than $2|\nu(\mcAf)|$, so the maximum number of extremal points is $n^{2|\nu(\mcAf)|}$ since there are $n$ possible choices for each coordinate (or IFS map) up to that length. \sqr

\section{Methods for the Convex Hull} \label{s04}

We proceed to the main results of this paper in determining the exact convex hull of polyfractals. ``Exactness'' is emphasized throughout to make it clear that the resulting extrema are the actual explicit extremal points, and not merely approximate or theoretical, as often is the case in the literature. The introduced ``Armadillo Method'' can be carried out in practice, as shown in the examples of Section \ref{s0504}.

As reasoned earlier, the case of irrational fractals seems hopelessly difficult, so we restrict our attention to fractals of unity, which by continuity in parameters approximate the irrational case (Theorem \ref{s030301}).

\subsection{A Method for Equiangular Fractals of Unity} \label{s0401}

Corollary \ref{s030205} and Theorems \ref{s030106}, \ref{s030204} hinted at the relevance of focal periodic extrema for the convex hull, but their fundamental role will only become clear in this section. We derive the convex hull of perhaps the simplest type of polyfractals, relevant to self-affine fractals.

\begin{prp} \label{s040101}
An $F=\langle T_1,\ldots,T_n\rangle,\ \vt_k=2\pi N/M\ (k\in\mcN),\ M\in\N,\ N\in\{0,\ldots,M-1\}$ equiangular fractal of unity can be generated as a Sierpi\'{n}ski fractal $F=\langle T_x: |x|=M'\rangle$ where $M'=M/\gcd(N,M)$, and furthermore $\Ext(F)=\Ext(p_x: |x|=M')$.\mydag
\end{prp}

The above implies that having equal IFS rotation angles $2\pi N/M$, it is sufficient to generate all $M'$-th level periodic points (computed using Corollary \ref{s020303}), as their convex hull will be that of the fractal. This is a simplified version of the method presented by K{\i}rat and Ko\c{c}yi\u{g}it \cite{ba00109} for the special case of planar equiangular IFS fractals.

To generalize the above to higher dimensions, define a focal address $x\in\mcAf$ in terms of matrix IFS factors $\cf_k\in\R^{d\times d}$ as one for which $\cf_x$ is a scalar multiple of the identity matrix. Then the part of Theorem \ref{s030106} guaranteeing $\Conv(F)=\Conv(\mcP)$ for Sierpi\'{n}ski fractals generalizes accordingly, to support the above proof.

\subsection{A General Method for Fractals of Unity} \label{s0402}

A general method for finding the convex hull of fractals of unity $F=\langle\mcT\rangle,\ \mcT=\{T_1,\ldots,T_n\}$ is presented, which have a finite number of extrema according to Theorem \ref{s030212}, though the method also works for rational fractals in general. The presented method is a reworked simplification of the one by K{\i}rat and Ko\c{c}yi\u{g}it \cite{ba00109} for the case of fractals of unity based on Section \ref{s0302}, the significant improvement being that termination is hereby guaranteed by Theorem \ref{s030211}.

Heavily utilizing the Extremal Theorem \ref{s040404} (ET), we take an ``inside-out blow-up approach'', meaning we attempt to find the convex hull by address generation. Proving rather inefficient, we take an ``outside-in'' approach with the Armadillo Method in Section \ref{s0404}, using the linear optimization algorithm of Appendix \ref{s0403}, making the two approaches philosophical antitheses.

According to Theorem \ref{s030211}, the extremal points of a fractal of unity have irreducible forms no longer than $2|\nu(\mcAf)|$ implying
\[ \Ext(F) = \Ext(T_b(p_x):\ |bx|\leq 2|\nu(\mcAf)|,\ \nu(x)=0,\ x\neq 0) \]
so the challenge becomes to compute the latter set efficiently, considering that its cardinality could reach $n^{2|\nu(\mcAf)|}$ making its computation potentially exponential in runtime. So we take a shortcut via irreducibility and a containment property.

\begin{dfn} \label{s040201}
An address $a\in\mcAf$ is \textbf{blowable} if $\exists b,x:\ bx<a, \nu(x)=0, x\neq 0$ denoted as $\exists\blo(a)$, and taking the shortest such $bx$ denote its \textbf{blow-up} as $\blo(a)\libeq T_b(p_x)$ (note that $(b,x)\in\Frm(\blo(a))$ is then irreducible). Denote the \textbf{blow-up of a set} $\mcA_0\subset\mcAf$ as $\blo(\mcA_0)\libeq \{\blo(a):\ a\in\mcA_0, \exists\blo(a)\}$. A set $\mcA_0\subset\mcAf$ is blowable if all of its elements are, denoted as $\exists\blo(\mcA_0)$. Define the set of \textbf{eventually focal points of a level} $L\in\N$ as the blow-up of all blowable $L$-long addresses $\EFoc_L(F)\libeq \{\blo(a):\ |a|=L,\ \exists\blo(a)\}$. (Note the property $\EFoc_l(F)\subset\EFoc_L(F)\ \forall l\leq L$.)
\end{dfn}

The set $\EFoc_L(F)$ can be easily generated via an algorithm which examines longer-and-longer addresses for eventual focality, and if a $(b,x)$ form is found (meaning $\nu(b)=\nu(bx),\ x\neq 0$) it will remain irreducible for higher iteration levels (see Section \ref{s0302}) so no further addresses $c>bx$ need to be examined (corresponding to sub-fractals of $T_bT_x(F)$), according to the ET. Iterating up to level $2|\nu(\mcAf)|$, the truncations of the extremal addresses necessarily emerge, due to Theorems \ref{s030211} and \ref{s030204}.

According to the Containment Lemma \ref{s020301}, if the containment property $\mH(S)\subset\Conv(S)$ holds for some compact $S\subset\C$, then $F=\langle\mcT\rangle\subset\Conv(S)$ implying that $\Conv(F)\subset\Conv(S)$. Furthermore if $S\subset F$ also holds, then $\Conv(S)\subset\Conv(F)$ so we can conclude that $\Conv(F)=\Conv(S)$ or equivalently that $\Ext(F)=\Ext(S)$. So we attempt to find such a set $S$ among the sets $\EFoc_L(F)\subset F$ for increasing iteration levels $L\in\N$. Theorem \ref{s030211} guarantees that such a level $L_*\leq 2|\nu(\mcAf)|$ will be found, implying that the method below is guaranteed to terminate with an output satisfying $\Ext(F)=\Ext(\EFoc_{L_*}(F))$.

\begin{mth} \label{s040202}\
A method for finding the exact convex hull of IFS fractals of unity.\myddag
\begin{enumerate}
\item Let the initial level be $L\libeq 1$.
\item Compute the set $\EFoc_L(F)$ and test whether the containment property $\mH(\EFoc_L(F))\subset\Conv(\EFoc_L(F))$ holds. If it does, then let $L_*\libeq L$ and go to Step 3. If it does not, then increase the level $L$ by one and repeat Step 2.
\item Conclude that $\Ext(F)=\Ext(\EFoc_{L_*}(F))$.
\end{enumerate}
\end{mth}

Compared to the method of K{\i}rat and Ko\c{c}yi\u{g}it \cite{ba00109}, the main difference in our method for IFS fractals of unity, is the increased efficiency due to the introduced idea of ``focality'' and the ET. Their method blows up addresses of the form $bx^lc,\ l\in\N$ (not considering the focality of $x$), while here a $bx,\ \nu(x)=0,\ x\neq 0$ form is sufficient for blow-up due to the ET. Computationally this is a significant improvement, since searching for a repeating part $x$ in some $a=bx^lc,\ |a|=L$ is more costly than reaching a point when $\nu(b)=\nu(a),\ |a|=L$ for some $b\lneq a$ in the address generation $n$-ary tree. This also excludes redundant eventually periodic addresses where the periodic part is not focal, in light of Theorem \ref{s030204}. The above change guarantees the termination of this method at some $L_*\leq 2|\nu(\mcAf)|$ due to Theorem \ref{s030211}, which was not guaranteed by their method.

Note that instead of the containment property of Step 2, their method essentially checks whether the condition $(\mH(E_L)\setminus E_L)\cup \Ext(F_L) \subset \Conv(E_L)$ holds, where $F_L\libeq \mH^L(\{p\}),\ p\in\mcP,\ E_L\libeq \EFoc_L(F)$. Here $E_L\subset \Conv(E_L)$ necessarily, so this condition can be simplified to $\mH(E_L)\cup \Ext(F_L)\subset\Conv(E_L)$. It is even less restrictive and sufficient to just require the containment property $\mH(E_L)\subset\Conv(E_L)$ as explained prior to the above method, since $S=E_L\subset F$ is compact. This is the containment checked in Step 2, eliminating the redundancy of finding $\Ext(F_L)$ and further set operations.

Furthermore, Method \ref{s040202} readily generalizes to higher dimensions as well, according to the remarks made after Proposition \ref{s040101} about defining focality with matrix IFS factors. Though our method may be appealing for its generality, it is inefficient for various reasons.

For instance, the computational cost of evaluating the sets $\EFoc_L(F)$ can be high in Step 2, potentially carried out for each $L\leq 2|\nu(\mcAf)|$. Though the property $\EFoc_l(F)\subset\EFoc_L(F)\ (l\leq L)$ and the remarks after Definition \ref{s040201} can increase the efficiency of evaluation cumulatively.

On the other hand, checking whether $\mH(\EFoc_L(F))\subset\Conv(\EFoc_L(F))$ holds can carry a significant risk of computational error. If the cardinality $|\Ext(F)|$ is large, then a program may arrive at the erroneous conclusion that the containment property is true prematurely, not at the actual sought level $L_*$ since the extrema tend to cluster, as noted in the remarks after Theorem \ref{s030110}.

Figures \ref{s050407} and \ref{s050410} generated from a seed, show why an ``inside-out'' address generation approach (Proposition \ref{s020204}) cannot possibly tackle the convex hull problem robustly. For certain IFS parameters, the extremal period length $|x|$ can be large (in these figures $262$ and $41$ respectively), potentially causing an exponential runtime for the iterative address generation of $\EFoc_L(F)$ (see Theorems \ref{s030211}, \ref{s030212}).

As projected earlier, due to the above issues with this ``inside-out'' method, an alternative ``outside-in'' approach is presented in Section \ref{s0404}.

\subsection{The Armadillo Method for Regular Fractals} \label{s0404}

The method to be introduced utilizes certain ideal target directions -- called ``regular directions'' -- to locate one focal extremal point, which is exploited via its cycle to generate a set of neighboring extrema. These extrema when iterated by each IFS map alone, will map out the entire convex hull. This method works if such a target direction exists, so in the next section we discuss a general heuristic candidate for bifractals, called the ``principal direction''. We continue to discuss fractals of unity $F=\langle\mcT\rangle,\ \mcT=\{T_1,\ldots,T_n\}$.

\begin{dfn} \label{s040401}
We say that two distinct extremal points are \textbf{neighbors} if the fractal is a subset of a closed half-space determined by the line connecting them. The \textbf{right/left neighbor} of an extremal point is the neighbor which comes counter/clockwise around the vertices of the convex hull. A set of extrema is \textbf{consecutive} if each element has a neighbor in the set -- all except two, have two -- and we call such a set a \textbf{plate}. A focal address $x\in\mcAf,\ \nu(x)=0$ is \textbf{consecutive} if $p_x\in\Ext(F)$ and $\Cyc(x)\subset\Ext(F)$ is consecutive around $\Conv(F)$; furthermore we say that $p_x$ is a \textbf{consecutive focal point}.
\end{dfn}

\begin{dfn} \label{s040405}
A target direction is \textbf{regular} if its maximizing extremal point is unique, strictly focal, and consecutive. A polyfractal is a \textbf{regular fractal} if it has a regular target direction, meaning if it has a strictly focal extremal point whose cycle is consecutive.
\end{dfn}

The regularity of a target direction $\tau\in\C$ can be verified by running Algorithm \ref{s040307}, which results in the maximizing irreducible form(s). If $(b,x),\ b\neq 0$ is one of these forms, then by the ET $p_x\in\Ext(F)$ is the maximizer of $\tau'\libeq \me^{-\nu(b)\vi}\tau$, so we can still check whether $x$ is strictly focal and $\Cyc(x)$ is consecutive, and we may conclude that $\tau'$ is regular instead.

Consecutiveness of the cycle can be verified as follows. To test if $e_{1,2}\in\Cyc(x)\subset\Ext(F)$ are neighbors, take an outward normal $\tau_0\in\C$ to the line connecting them, and run Algorithm \ref{s040307} to see what irreducible form(s) maximize(s) this target. If those of $e_1$ and $e_2$ do, then they must be neighbors. But how does one determine the order of the elements of $\Cyc(x)\subset\Ext(F)$ on the boundary of the convex hull? This would be a prerequisite for testing the cycle for consecutiveness. Simply connect each element $e\in\Cyc(x)$ to some fixed interior point of the convex hull $f\in\mint\ \Conv(F)$, and the values $\arg(e-f)$ imply the sought order.

\begin{thm} \label{s040408}
Let $p_x\in\Ext(F)\cap\Foc(F)$ be the maximizer of a regular target direction. Then the iterates of the plate $\Cyc(x)$ by each $T_k\in\mcT$ generate all the extrema of $F$, meaning
\[ \Conv(F) = \Conv\ \bigcup_{k=1}^n \bigcup_{l=0}^{\left\lceil\frac{\pi}{|\vt_k|}\right\rceil} T_k^l\left(\Cyc(x)\right). \]
\end{thm}
\noindent
\prf
Due to regularity we have that $\Cyc(x)\subset\Ext(F)$ are consecutive extrema. By the strict focality of $x$ we have that each $k\in\mcN=\{1,\ldots,n\}$ occurs in $x$ so $\forall k\in\mcN\ \exists e_{1,2}\in\Cyc(x):\ e_1=T_k(e_2)$ therefore the angle with vertex at $p_k$ and spanned by $\Cyc(x)$ is at least $|\vt_k|$ for each $k\in\mcN$. So the iterated plate $T_k^l(\Cyc(x)), l\in\N$ produces more (potentially overlapping) plates of extrema. Since each $T_k$ traces a logarithmic spiral trajectory, we only need to iterate up to a $\left\lceil\pi/|\vt_k|\right\rceil$ number of times for each $k\in\mcN$. This gives the convex hull identity in the theorem. The identity also readily implies that $F$ has a finite number of extrema, since on the left of the identity we have a finite union of finite sets. \sqr

\begin{mth} \label{s040410} \textup{(Armadillo Method)}\
Determines the extremal points of a fractal of unity.\myddag
\begin{enumerate}
\item Find a reasonable candidate direction $\tau\in\C\setminus\{0\}$ which is potentially regular.
\item Run Algorithm \ref{s040307} (LOAF) and see if it gives a unique strictly focal maximizing form $(0,x)$. If it does, then proceed to Step 3. Otherwise go to Step 1 and try a different direction.
\item Evaluate the cycle $\Cyc(x)\subset\Ext(F)$ and deduce the order of its elements on the boundary of the convex hull, as discussed above.
\item Determine if $\Cyc(x)$ is a plate, by connecting its elements with lines (ordered in the previous step), and running the LOAF. If the algorithm returns the two connected cycle extrema, then they are neighbors on the convex hull.
\item Iterate the cycle $\Cyc(x)$ by each IFS map according to Theorem \ref{s040408} to find the rest of the extremal points.
\end{enumerate}
\end{mth}

In the last step, the iterates of the plate $\Cyc(x)$ can overlap -- an easily resolvable redundancy -- similarly to the plates of the armored mammal armadillo. The above method is potentially more efficient and robust than the general method of Section \ref{s0402}. A candidate target direction can be readily verified for regularity with the LOAF. While the general method blows up all finite addresses up to a level (inefficient, since it may run up to an address length of $2|\nu(\mcAf)|$ making it exponential) and checks a containment property for the blow-ups (may not be robust), the LOAF efficiently eliminates most subfractals recursively.

An approximative method for finding a regular direction, could be to find one for a ``simpler'' fractal. Approximate the IFS rotation angles $\vt_k/2\pi$ with low-numerator fractions, and find a regular direction for this simpler IFS (Section \ref{s0502} reasons why numerators are relevant). By Theorem \ref{s030301}, the extrema of the simplified fractal approximate that of the original, so if the approximation is not too crude, then the fractals should share this regular direction.

\section{The Principal Direction} \label{s05}

In this section, we restrict our attention to ``C-IFS'' fractals of unity in our quest to find a reasonable candidate target direction -- called the ``principal direction'' -- the regularity of which can then be verified via the methods discussed in the previous section. Nevertheless a detailed heuristic reasoning is given for its probable regularity, while a solid proof may not even be possible, since there can be rare counterexamples.

\subsection{The Normal Form of Bifractals} \label{s0501}

To simplify a discussion, it is often preferable to transform a bifractal $F=\langle T_1,T_2\rangle$ to ``normal form'', i.e. to normalize the fixed points $p_{1,2}\in\mcP$ of the IFS to $0$ and $1$ (it is up to our preference which fixed point we normalize to $0$). This can be accomplished with the affine similarity transform $N(z)\libeq (z-p_1)/(p_2-p_1)$.

\begin{prp} \label{s050101}
With the above $N$ map, we have for $T_k'\libeq N\co T_k\co N^{-1},\ k=1,2$ that $T_1'(z)= \cf_1 z,\ T_2'(z)= 1+\cf_2(z-1)$ and $N(F)=\langle T_1', T_2'\rangle$ which we call the \textbf{normal form} of $F$. Furthermore for any $x\in\mcAf$ using $T_x'= N T_x N^{-1}$ we can express
\[ p_x = \frac{T_x(0)}{1-\cf_x} = p_1 + \frac{T_x'(0)}{1-\cf_x}(p_2-p_1).\mydag \]
\end{prp}

Essentially, the above proposition implies that we can scale down the IFS fractal to normal form, examine its geometrical properties (such as its convex hull), and then transform back the results, since due to the Affine Lemma \ref{s020305} $F=N^{-1}\langle N\co T_1\co N^{-1}, N\co T_2\co N^{-1}\rangle$. Therefore the factors $\cf_{1,2}$ entirely characterize the geometry of a bifractal, since they remain unchanged under normalization. Thus we restrict the discussion to normalized bifractals and drop the apostrophes, so the IFS maps become $T_1(z)=\cf_1 z$ and $T_2(z)=1+\cf_2(z-1)$.

\begin{prp} \label{s050102}
All finite compositions of $T_{1,2}$ ending in $T_2$ with seed $p_1=0$ take the form
\[ T_1^{n_0}T_2T_1^{n_1}T_2\ldots T_1^{n_L}T_2(0) = (1-\cf_2)\sum_{j=0}^L \cf_1^{s_j}\cf_2^j\ \ \mathrm{where}\ \ L\in\N,\ n_j\geq 0,\ s_j\libeq \sum_{k=0}^j n_k.\mydag \]
\end{prp}

\subsection{A Candidate Direction for C-IFS Fractals} \label{s0502}

\begin{dfn} \label{s050201}
We say that a bifractal is a \textbf{C-IFS fractal} (after the L\'{e}vy C curve), if its rotation angles satisfy $\vt_1\in(-\pi,0),\ \vt_2\in(0,\pi),\ |\vt_1|\leq\vt_2$.
\end{dfn}

\begin{heu} \label{s050202}
Let $F=\langle T_1,T_2\rangle$ be a normalized C-IFS fractal of unity with
\[ \vt_1=-\frac{2\pi P}{M},\ \vt_2=\frac{2\pi Q}{M}\ \ \mathrm{where}\ \ P,Q,M\in\N,\ 0<P\leq Q<\frac{M}{2}. \]
Then the target direction $\tau_*\libeq \vi(1-\cf_2)\Log\ \cf_1$ called the \textbf{principal direction}, is likely to be regular, and its maximizer is likely to be the strictly focal periodic point $p_x$ where
\[ T_x=T_2T_1^{n_1}\ldots T_2T_1^{n_J},\ \ |x|=\frac{P+Q}{\gcd(P,Q)} \]
\[ n_j=s_j-s_{j-1},\ j=1,\ldots,J\libeq \frac{P}{\gcd(P,Q)},\ n_0=s_0=0 \]
\[ s_j\libeq \argmax\left(\vl_1^s\cos(\vt_1 s+\vt_2 j+\alpha):\ s=\left\lfloor\frac{Qj}{P}\right\rfloor,\ \left\lceil\frac{Qj}{P}\right\rceil\right),\ \ \alpha\libeq \arctan\left(\frac{\ln \vl_1}{\vt_1}\right). \]
Let a maximizer with respect to the principal direction be called a \textbf{principal extremal point}, and denote it as $e_*$ if it is unique.
\end{heu}
\noindent
\rea
First of all, we calculate the Euclidean inner product of the principal direction $\tau_*$ with a general normalized fractal point given by Proposition \ref{s050102}, and connect it to $s_j$. Our calculation is aided by the identity $\forall u,v\in\C: \langle u,v\rangle=\mRe(u\bar{v})$. Denoting $t_j(s)\libeq \vl_1^s\cos(\vt_1 s+\vt_2 j+\alpha)$ we have
\[ \left\langle\vi(1-\cf_2)\Log\ \cf_1,\ (1-\cf_2)\sum_{j=0}^L \cf_1^{s_j}\cf_2^j\right\rangle = |1-\cf_2|^2\sum_{j=0}^L \left\langle\vi\ \Log\ \cf_1,\ \cf_1^{s_j}\cf_2^j\right\rangle = \]
\[ = |1-\cf_2|^2|\Log\ \cf_1| \sum_{j=0}^L \vl_1^{s_j}\vl_2^j\cos \angle\left(\vi\ \Log\ \cf_1,\ \cf_1^{s_j}\cf_2^j\right) = |1-\cf_2|^2|\Log\ \cf_1| \sum_{j=0}^L \vl_2^j t_j(s_j) \]
\[ \mathrm{since}\ \ \angle\left(\vi\ \Log\ \cf_1,\ \cf_1^{s_j}\cf_2^j\right) \equiv \Arg(\cf_1^{s_j}\cf_2^j) - \Arg(-\vt_1+\vi\ln \vl_1) \equiv \]
\[ \equiv (\vt_1s_j+\vt_2j)-\arctan\left(-\frac{\ln \vl_1}{\vt_1}\right) \equiv \vt_1s+\vt_2j+\alpha\ (\mmod\ 2\pi). \]
So we deduced that in order to maximize $z\mapsto\langle\tau_*,z\rangle$ over $F=\langle T_1,T_2\rangle$ we need to maximize the sum $\sum_{j=0}^L \vl_2^j t_j(s_j)$ over all $L\in\N$ and $0\leq s_j\leq s_{j+1}\ (j=0,1,2,\ldots)$.\\
We temporarily drop the latter optimization constraint, and just consider each $s\mapsto t_j(s)$ function individually, and see where its local extrema are. Taking the derivative
\[ t_j'(s) = (\ln \vl_1)\vl_1^s\cos(\vt_1s+\vt_2j+\va)-\vt_1\vl_1^s\sin(\vt_1s+\vt_2j+\va) = \]
\[ = (\ln \vl_1)\vl_1^s\left[\cos(\vt_1s+\vt_2j)\cos(\va)-\sin(\vt_1s+\vt_2j)\sin(\va)\right] - \]
\[ -\ \vt_1\vl_1^s \left[\sin(\vt_1s+\vt_2j)\cos(\va)+\cos(\vt_1s+\vt_2j)\sin(\va)\right] = \]
\[ = \vl_1^s\left(\cos(\vt_1s+\vt_2j)\left[(\ln \vl_1)\cos \va - \vt_1\sin \va\right] - \sin(\vt_1s+\vt_2j)\left[(\ln \vl_1)\sin \va + \vt_1\cos \va\right]\right) = \]
\[ = -\vl_1^s\sin(\vt_1s+\vt_2j)\left[(\ln \vl_1)\sin \va + \vt_1\cos \va\right] = \]
\[ = |\vt_1|(\cos \va)(1+(\tan \va)^2) \sin(\vt_1s+\vt_2j) \]
\[ \mathrm{since}\ \ (\ln \vl_1)\cos \va - \vt_1\sin \va = (\cos \va)\vt_1(\tan \va - \tan \va) = 0 \]
\[ \mathrm{and}\ \ (\ln \vl_1)\sin \va + \vt_1\cos \va = \vt_1(\cos \va)((\tan \va)^2 +1)\neq 0. \]
Here $\cos \va>0$ by $\alpha=\arctan\left(\frac{\ln \vl_1}{\vt_1}\right)\in\left(0,\frac{\pi}{2}\right)$ due to $\frac{\ln \vl_1}{\vt_1}>0$ so remarkably
\[ \msgn\ t_j'(s) = \msgn\ \sin(\vt_1s + \vt_2j) = -\msgn\ \sin(|\vt_1|s-\vt_2j) \]
which ultimately followed from our special choice of $\tau_*$. Since the $\vt_1<0$ factor flips the oscillation of $\sin$, the local maxima occur at
\[ \vt_1s+\vt_2j=2\pi k,\ k\in\Z\ \ \rar\ \ s\in\left\{-\frac{\vt_2 j}{\vt_1}+\frac{2\pi}{\vt_1} k:\ k\in\Z\right\} = \left\{ \frac{Qj}{P}+\frac{M}{P} k:\ k\in\Z\right\}. \]
Denote by $s_j$ the heuristically most reasonable sequence $j=0,1,2,\ldots$ for maximizing $f\mapsto\langle\tau_*,f\rangle,\ f\in F$. For $j=0$ the likely solution is $s_0=0$ with $k=0$, since for $k<0$ the number $\frac{Qj}{P}+\frac{M}{P}k$ is negative, and for $k\geq 1$ the maximal values of $t_j$ decrease, since $s\mapsto t_0(s)$ oscillates between the decreasing exponential curves $s\mapsto\pm\vl_1^s$. Therefore we conclude heuristically that $s_0=0$.\\
We further reason that $s_j\in\left\{\left\lfloor\frac{Qj}{P}\right\rfloor, \left\lceil\frac{Qj}{P}\right\rceil\right\}$ for each $j\in\N$. Since $\frac{Qj}{P}$ occur with a spacing of $\frac{Q}{P}<\frac{M}{2P}$ while the local maxima of $t_j$ occur with a spacing of $\frac{M}{P}$ then assuming heuristically the somewhat stricter constraint $s_j<s_{j+1},\ j\in\N$ we necessarily have that $s_j$ occurs in an integer near $\frac{Qj}{P}$.\\
We may thus conclude the heuristic statement that the maximizing address with respect to $\tau_*$ is likely to have the collected exponents $n_j=s_j-s_{j-1},\ j\in\N$ where
\[ s_j=\argmax\left(t_j(s):\ s=\left\lfloor\frac{Qj}{P}\right\rfloor,\ \left\lceil\frac{Qj}{P}\right\rceil\right). \]
We show that this address is periodic according to the map
\[ T_x=T_2T_1^{n_1}\ldots T_2T_1^{n_J},\ j=1,\ldots,J=\frac{P}{\gcd(P,Q)} \]
or equivalently that $n_j$ is periodic by $J$. If we showed that $s_{j+J}=s_j+s_J$ then it would imply periodicity $n_{j+J}=n_j$ since
\[ n_{j+J} = s_{j+J}-s_{j+J-1} = \left(s_j+s_J\right)-\left(s_{j-1}+s_J\right) = s_j-s_{j-1} = n_j. \]
First of all, $s_J=\frac{QJ}{P}$ necessarily, so $\vt_1 s_J+\vt_2 J=0$ implying
\[ t_{j+J}(s+s_J)= \vl_1^{s+s_J}\cos(\vt_1(s+s_J)+\vt_2(j+J)+\va)=\vl_1^{s_J}t_j(s). \]
Since the optimum of $t_j$ with respect to the global target $\tau_*$ was deduced to likely be $s_j\in\left\{\left\lfloor\frac{Qj}{P}\right\rfloor, \left\lceil\frac{Qj}{P}\right\rceil\right\}$ and due to
\[ s_{j+J}\in\left\{\left\lfloor\frac{Q(j+J)}{P}\right\rfloor,\ \left\lceil\frac{Q(j+J)}{P}\right\rceil\right\} = s_J + \left\{\left\lfloor\frac{Qj}{P}\right\rfloor,\ \left\lceil\frac{Qj}{P}\right\rceil\right\} \]
we necessarily have $s_{j+J}=s_j+s_J$ implying the periodicity of $n_j$ by $J$. Thus the maximizer of $\tau_*$ is heuristically likely to be $p_x$ which is strictly focal since $\nu(x)=\vt_1 s_J+\vt_2 J=0$ and $T_x$ contains both $T_1$ and $T_2$. Using Algorithm \ref{s040307}, we can verify if this $p_x$ is indeed the maximizer of $\tau_*$, in which case then necessarily $\Cyc(x)\subset\Ext(F)$. With this algorithm, we can also verify the consecutiveness of the cycle (see the remarks after Definition \ref{s040405}), which was typically the case in our experiments (i.e. $\tau_*$ was a regular target direction). \sqr

Note that due to $P,Q\leq M/2$ we can estimate
\[ |x| = \frac{P+Q}{\gcd(P,Q)} \leq \frac{M}{\gcd(P,Q,M)}\frac{\gcd(P,Q,M)}{\gcd(P,Q)}\leq \frac{M}{\gcd(P,Q,M)} = |\nu(\mcAf)| \]
according to Lemma \ref{s030202}, so Theorem \ref{s030211} suggests that $p_x$ could be in irreducible form.

\subsection{The Armadillo Method Adopted for C-IFS Fractals} \label{s0503}

\begin{mth}  \label{s050301} \textup{(Exact)}\
An exact method for finding the convex hull of C-IFS fractals of unity, based on the heuristic prediction that the principal direction is likely to be regular.\myddag
\begin{enumerate}
\item Normalize the IFS maps to $T_1(z)=\cf_1 z,\ T_2(z)=1+\cf_2(z-1)$.
\item Calculate the principal direction $\tau_*=\vi(1-\cf_2)\Log\ \cf_1$.
\item Run Algorithm \ref{s040307} (LOAF) and see if the maximizer of $\tau_*$ is a unique $p_x$ strictly focal point. If it is then proceed to Step 4, otherwise try a different candidate direction.
\item Next deduce the order of $e\in\Cyc(x)\subset\Ext(F)$ on the convex hull from the values $\arg(e-f)$ where $f\in\mint\ \Conv(F)$ is fixed. Verify if the cycle is consecutive by connecting each cycle point with its likely neighbor, and check using the LOAF in the normal direction to the connecting line, whether the maximizers are these two extrema. If the cycle turns out not to be consecutive, then try a different candidate direction in Step 3.
\item Iterate the cycle $\Cyc(x)$ by each IFS map to find the rest of the extrema (Theorem \ref{s040408}).
\item Lastly, map the determined extrema back by the inverse of the normalizing map $N^{-1}(z)=p_1+z(p_2-p_1)$ to get the extrema of the original fractal.
\end{enumerate}
If all steps can be carried out, then it results in the convex hull of a C-IFS fractal of unity.
\end{mth}

\begin{mth}  \label{s050302} \textup{(Heuristic)}\
A heuristic method based on Heuristic \ref{s050202} for finding the convex hull of C-IFS fractals of unity.
\begin{enumerate}
\item Normalize the IFS maps to $T_1(z)=\cf_1 z,\ T_2(z)=1+\cf_2(z-1)$.
\item Calculate the principal direction $\tau_*=\vi(1-\cf_2)\Log\ \cf_1$ and $J=P/\gcd(P,Q)$.
\item Calculate $s_j=\argmax\left(\vl_1^s\cos(\vt_1 s+\vt_2 j+\alpha):\ s=\left\lfloor\frac{Qj}{P}\right\rfloor, \left\lceil\frac{Qj}{P}\right\rceil\right),\ 1\leq j\leq J$.
\item Calculate $n_j=s_j-s_{j-1}$ and $x$ such that $T_x=T_2T_1^{n_1}\ldots T_2T_1^{n_J}$ and the resulting $p_x$ using Corollary \ref{s020303}.
\item Iterate the cycle $\Cyc(x)$ by each IFS map, according to the heuristic application of the Armadillo Method, to find the rest of the extrema.
\item Lastly, map the determined extrema back by the inverse of the normalizing map to get the extrema of the original fractal.
\end{enumerate}
This method results in the vertices of a polygon, which is heuristically predicted to be the convex hull of the C-IFS fractal.
\end{mth}

According to Theorem \ref{s030301}, for irrational bifractals one may resort to finding the extrema of increasingly better approximations by fractals of unity, varying the angular parameters. Considering Method \ref{s050302}, this means that the exponents $n_j$ will be generated ad infinitum.

\subsection{Examples} \label{s0504}

\begin{exa} \label{s050401} \textup{(L{\'e}vy C Curve \cite{ba00053, ba00054, ba00055, ba00056}, Figure \ref{s050404})}\\
Find the extrema of the C-IFS fractal of unity in normal form with IFS factors
\[ \cf_1=\frac{1}{\sqrt{2}}\exp\left(-\frac{\pi}{4}\vi\right),\ \ \cf_2=\frac{1}{\sqrt{2}}\exp\left(\frac{\pi}{4}\vi\right). \]
\end{exa}
\noindent
\sln
Applying Method \ref{s050301}, we maximize in the principal direction
\[\tau_*=\vi(1-\cf_2)\Log\ \cf_1\approx 0.2194-0.5660\vi \]
and using Algorithm \ref{s040307} we arrive at the unique strictly focal maximizing truncation $x=(2,1)$ resulting in the principal extremal point and its cycle
\[ e_*= p_x= \frac{T_x(0)}{1-\cf_x} = \frac{1-\cf_2}{1-\cf_1\cf_2} = 1-\vi \]
\[ \Cyc(x)=\{e_*, T_1(e_*)\}= \left\{\frac{1-\cf_2}{1-\cf_1\cf_2},\ \frac{\cf_1(1-\cf_2)}{1-\cf_1\cf_2}\right\}= \{1-\vi,\ -\vi\}. \]
Next deduce using Algorithm \ref{s040307} that the cycle is consecutive (meaning $e_*, T_1(e_*)$ are neighbors). Lastly, we iterate the cycle according to Theorem \ref{s040408} and get the rest of the extremal points $\Ext(F)=\{e_*,\ T_1(e_*),\ T_1^2(e_*),\ T_1^3(e_*),\ T_1^4(e_*),\ T_2(e_*),\ T_2^2(e_*),\ T_2^3(e_*)\}$. \sqr

\begin{exa} \label{s050402} \textup{(Twindragon / Davis--Knuth Dragon \cite{ba00057}, Figure \ref{s050405})}\\
Find the extrema of the C-IFS fractal of unity in normal form with IFS factors
\[ \cf_1=\frac{1}{\sqrt{2}}\exp\left(-\frac{\pi}{4}\vi\right),\ \ \cf_2=\frac{1}{\sqrt{2}}\exp\left(\frac{3\pi}{4}\vi\right). \]
\end{exa}
\noindent
\sln
Applying Method \ref{s050301}, we maximize in the principal direction
\[ \tau_*=\vi(1-\cf_2)\Log\ \cf_1\approx 1.0048-0.9126\vi \]
and using Algorithm \ref{s040307} we arrive at the unique strictly focal maximizing truncation $x=(2,1,1,1)$ or $T_x= T_2 T_1^3$ resulting in the principal extremal point and its cycle
\[ e_*= p_x= \frac{T_x(0)}{1-\cf_x} = \frac{1-\cf_2}{1-\cf_1^3\cf_2} = 2-\frac23\vi \]
\[ C_x\libeq \Cyc(x)= \{e_*,\ T_1(e_*),\ T_1^2(e_*),\ T_1^3(e_*)\} = \]
 \[ = \left\{\frac{1-\cf_2}{1-\cf_1^3\cf_2},\ \frac{\cf_1(1-\cf_2)}{1-\cf_1^3\cf_2},\ \frac{\cf_1^2(1-\cf_2)}{1-\cf_1^3\cf_2},\ \frac{\cf_1^3(1-\cf_2)}{1-\cf_1^3\cf_2}\right\} = \]
 \[ = \left\{2-\frac23\vi,\ \frac23-\frac43\vi,\ -\frac13-\vi,\ -\frac23-\frac13\vi\right\}. \]
Next deduce using Algorithm \ref{s040307} that the cycle is consecutive. Lastly, we iterate the cycle according to Theorem \ref{s040408} and get the rest of the extremal points $\Ext(F)\subset C_x\cup T_1(C_x)\cup T_2(C_x)$. \sqr

\begin{exa} \label{s050403} \textup{(Figure \ref{s050406})}\\
Find the extrema of the C-IFS fractal of unity in normal form with IFS factors
\[ \cf_1=0.65\exp\left(-\frac{2\pi}{6}\vi\right),\ \ \cf_2=0.65\exp\left(\frac{3\pi}{6}\vi\right). \]
\end{exa}
\noindent
\sln
Applying Method \ref{s050301}, we maximize in the principal direction
\[ \tau_*=\vi(1-\cf_2)\Log\ \cf_1\approx 0.7672-1.1115\vi \]
and using Algorithm \ref{s040307} we arrive at the unique strictly focal maximizing truncation $x=(2,1,1,2,1)$ or $T_x= T_2 T_1^2 T_2 T_1$ resulting in the principal extremal point and cycle
\[ e_*= p_x= \frac{T_x(0)}{1-\cf_x} = \frac{(1-\cf_2)(1+\cf_1^2\cf_2)}{1-\cf_1^3\cf_2^2}\approx 1.2993-1.0655\vi \]
\[ C_x\libeq \Cyc(x)= \{e_*,\ T_1(e_*),\ T_2 T_1(e_*),\ T_1 T_2 T_1(e_*),\ T_1^2 T_2 T_1(e_*)\}. \]
Next deduce using Algorithm \ref{s040307} that the cycle is consecutive. Then iterate the cycle according to Theorem \ref{s040408} and get the rest of the extremal points $\Ext(F)\subset C_x\cup T_1(C_x)\cup T_1^2(C_x)\cup T_2(C_x)\cup T_2^2(C_x)$.\ \sqr

In the figures that follow, the principal extremal point $e_*$ is marked by a blue dot, and the blue line is perpendicular to the principal direction. The fixed points $p_{1,2}$ are plotted with red $\odot$ and the iterates $T_1(p_2), T_2(p_1)$ with magenta. The cycle vertices are circled in red. Figures \ref{s050407}, \ref{s050410}, and \ref{s050411} clearly illustrate the utility of the Armadillo Method \ref{s040410}.

The IFS factors and the principal extrema for the rest of the figures:
\begin{itemize}
\item Figure \ref{s050407}:\tabto{2.05cm}$\cf_1=0.6177\exp\left(-\frac{99\pi}{180}\vi\right),\ \ \cf_2=0.8594\exp\left(\frac{163\pi}{180}\vi\right)$
    \vspace{0.1cm}
    \tabto{2cm}$T_x=T_2 (T_a^4T_b)^2 (T_a^5T_b)^4T_aT_2^{-1},\ \ e_*=p_x\approx 2.7130-0.5959\vi$
    \vspace{0.1cm}
    \tabto{2cm}with\ $T_a=T_1^2T_2T_1T_2T_1^2T_2,\ \ T_b=T_1T_2T_1^2T_2$\ and
    \vspace{0.1cm}
    \tabto{2cm}$|x|=\frac{99+163}{\gcd(99,163)}=1+2(4|a|+|b|)+4(5|a|+|b|)+|a|-1=29|a|+6|b|=262$.
\item Figure \ref{s050408}:\tabto{2.05cm}$\cf_1=0.6\exp\left(-\frac{8\pi}{12}\vi\right),\ \ \cf_2=0.8\exp\left(\frac{9\pi}{12}\vi\right)$
    \vspace{0.1cm}
    \tabto{2cm}$|x|=\frac{8+9}{\gcd(8,9)}=17,\ \ T_x=(T_2T_1)^4 T_1(T_2T_1)^4,\ \ e_*=p_x\approx 2.9439-0.5767\vi$.
\item Figure \ref{s050409}:\tabto{2.05cm}$\cf_1=0.5479\exp\left(-\frac{5\pi}{45}\vi\right),\ \ \cf_2=0.9427\exp\left(\frac{12\pi}{45}\vi\right)$
    \vspace{0.1cm}
    \tabto{2cm}$|x|=\frac{5+12}{\gcd(5,12)}=17,\ \ T_x= (T_2T_1^2 T_2T_1^3)^2 T_2T_1^2,\ \ e_*=p_x\approx 0.5203-0.9244\vi$.
\item Figure \ref{s050410}:\tabto{2.05cm}$\cf_1=0.9\exp\left(-\frac{6\pi}{45}\vi\right),\ \ \cf_2=0.5\exp\left(\frac{35\pi}{45}\vi\right)$
    \vspace{0.1cm}
    \tabto{2cm}$|x|=\frac{6+35}{\gcd(6,35)}=41,\ \ T_x= (T_2T_1^6)^3 T_2T_1^5(T_2T_1^6)^2,\ \ e_*=p_x\approx 1.8720-0.4808\vi$.
\item Figure \ref{s050411}:\tabto{2.05cm}$\cf_1=0.9421\exp\left(-\frac{2\pi}{180}\vi\right),\ \ \cf_2=0.9561\exp\left(\frac{17\pi}{180}\vi\right)$
    \vspace{0.1cm}
    \tabto{2cm}$|x|=\frac{2+17}{\gcd(2,17)}=19,\ \ T_x= T_2T_1^9 T_2T_1^8,\ \ e_*=p_x\approx 0.1958-0.6532\vi$.
\end{itemize}

\section{Concluding Remarks} \label{s06}

Methods have been introduced to determine the exact (not approximate) convex hull of certain classes of IFS fractals, and it has been reasoned how this classification prevents such efforts beyond the class of ``rational'' IFS fractals. A general though inefficient method was discussed for fractals of unity, and a more efficient method for regular fractals. In the latter ``Armadillo Method'', a certain periodic extremal point was exploited to generate the rest of the convex hull.

These methods were designed with potential generalization to 3D IFS fractals in mind -- as remarked after Proposition \ref{s040101} -- so that the convex hull of plants like L-system trees or the Romanesco Broccoli can be found. In fact with those remarks, Sections \ref{s03} and \ref{s04} can be generalized to higher dimensions. The primary aim of this paper was to find the most natural approach to the Convex Hull Problem in the plane.

\begin{appendices}

\section{Example Figures} \label{s0701}

\begin{figure}[H]
\centering
\includegraphics[width=360pt]{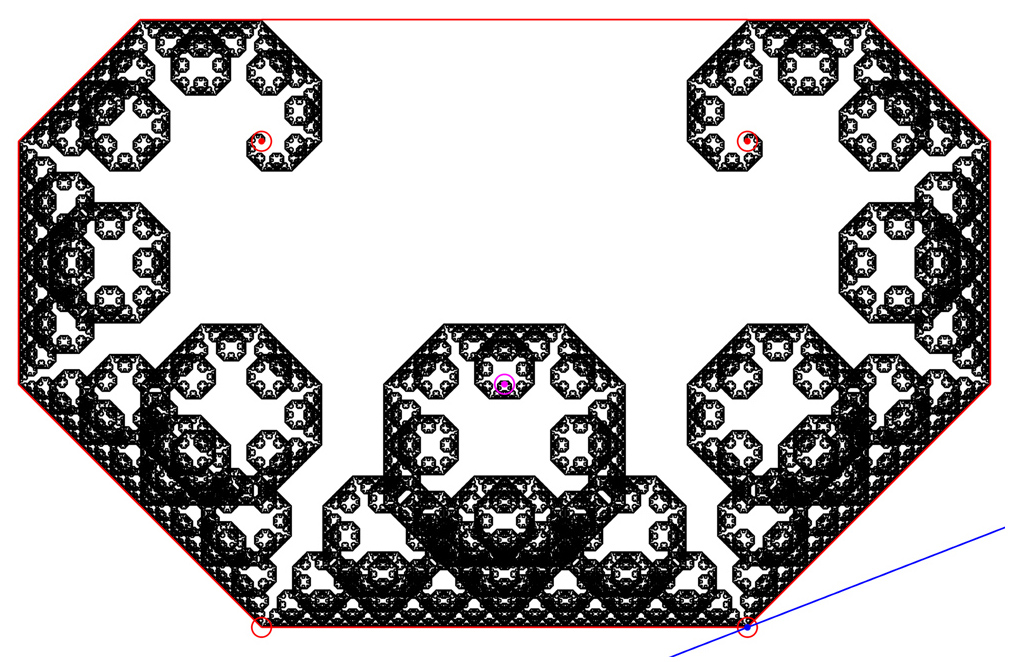}
\caption{L{\'e}vy C Curve, Example \ref{s050401}.}
\label{s050404}
\end{figure}
\begin{figure}[H]
\centering
\includegraphics[width=360pt]{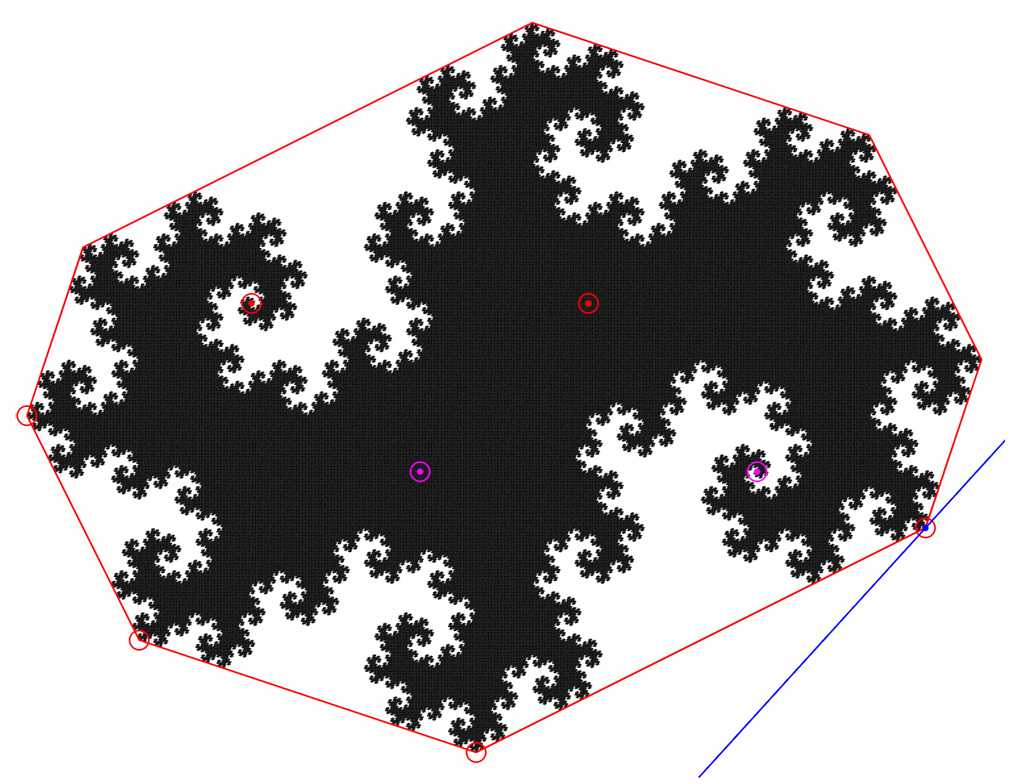}
\caption{Twindragon, Example \ref{s050402}.}
\label{s050405}
\end{figure}
\begin{figure}[H]
\centering
\includegraphics[width=370pt]{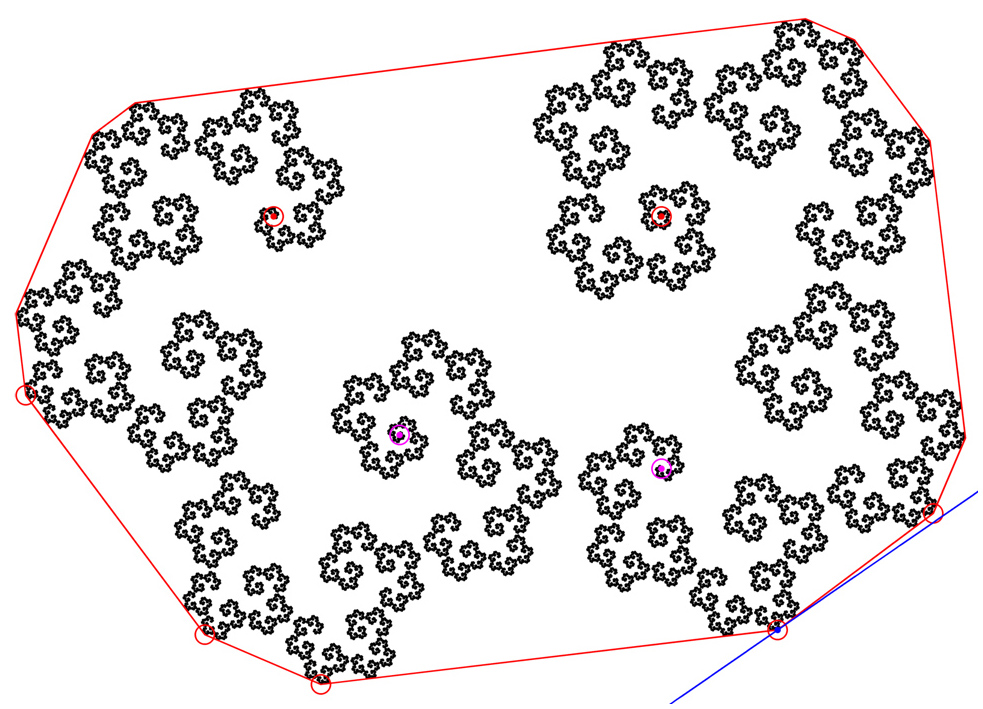}
\caption{Example \ref{s050403}.}
\label{s050406}
\end{figure}
\begin{figure}[H]
\centering
\includegraphics[width=370pt]{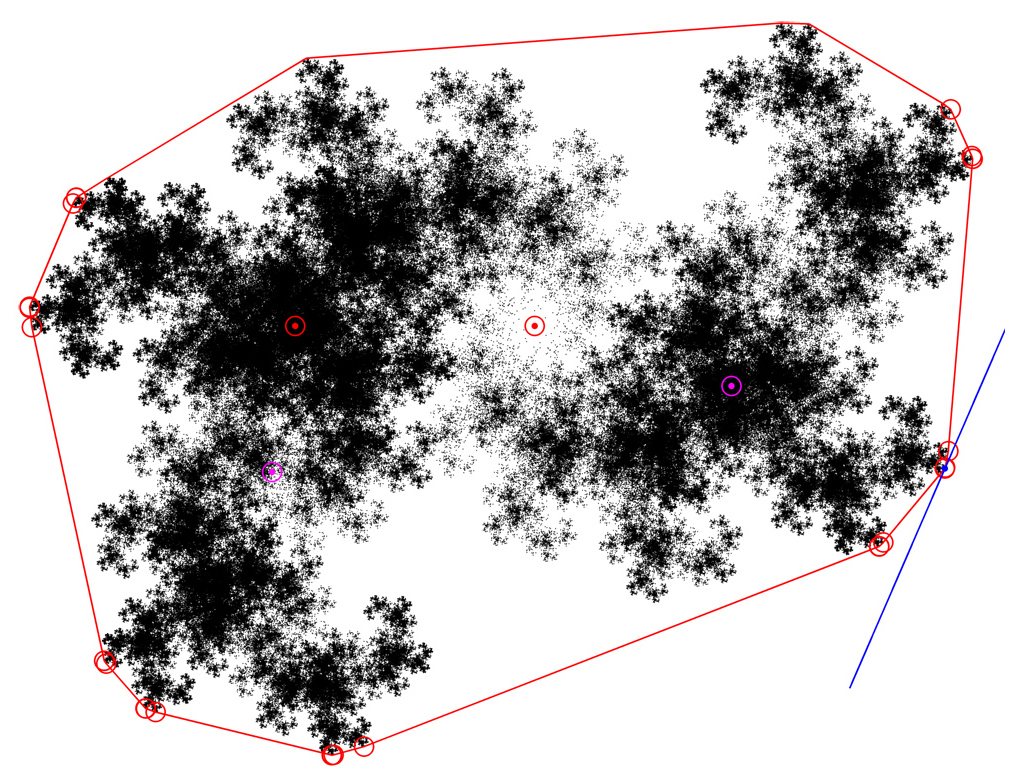}
\caption{Illustration of the subtlety of convex hull determination: $|\Ext(F)|\geq 262$.}
\label{s050407}
\end{figure}
\begin{figure}[H]
\centering
\includegraphics[width=350pt]{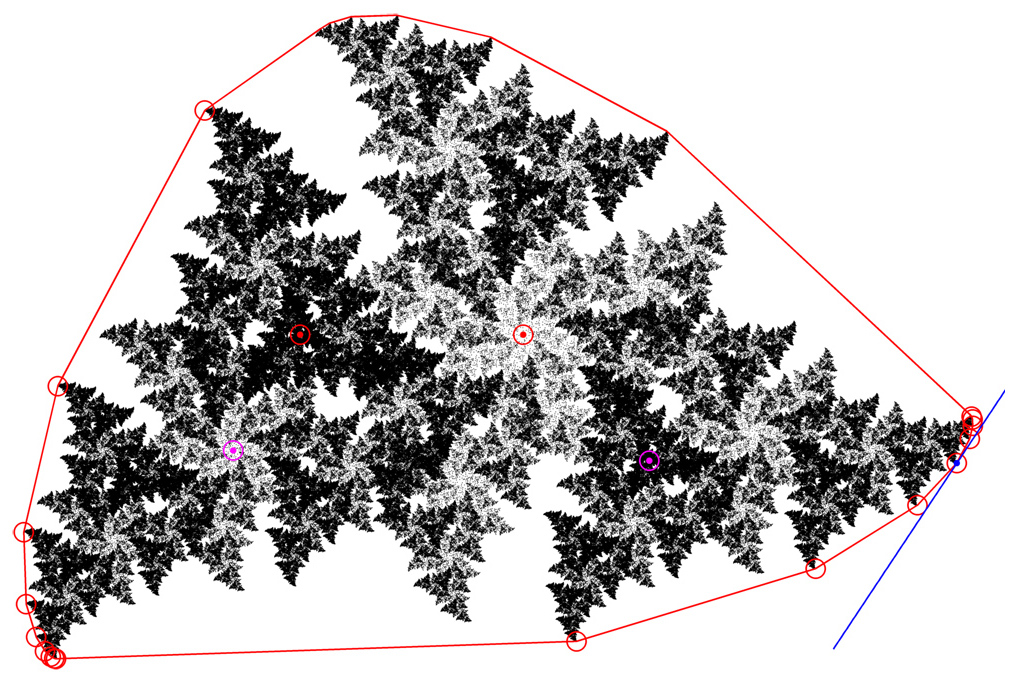}
\caption{A random C-IFS fractal.}
\label{s050408}
\end{figure}
\begin{figure}[H]
\centering
\includegraphics[width=350pt]{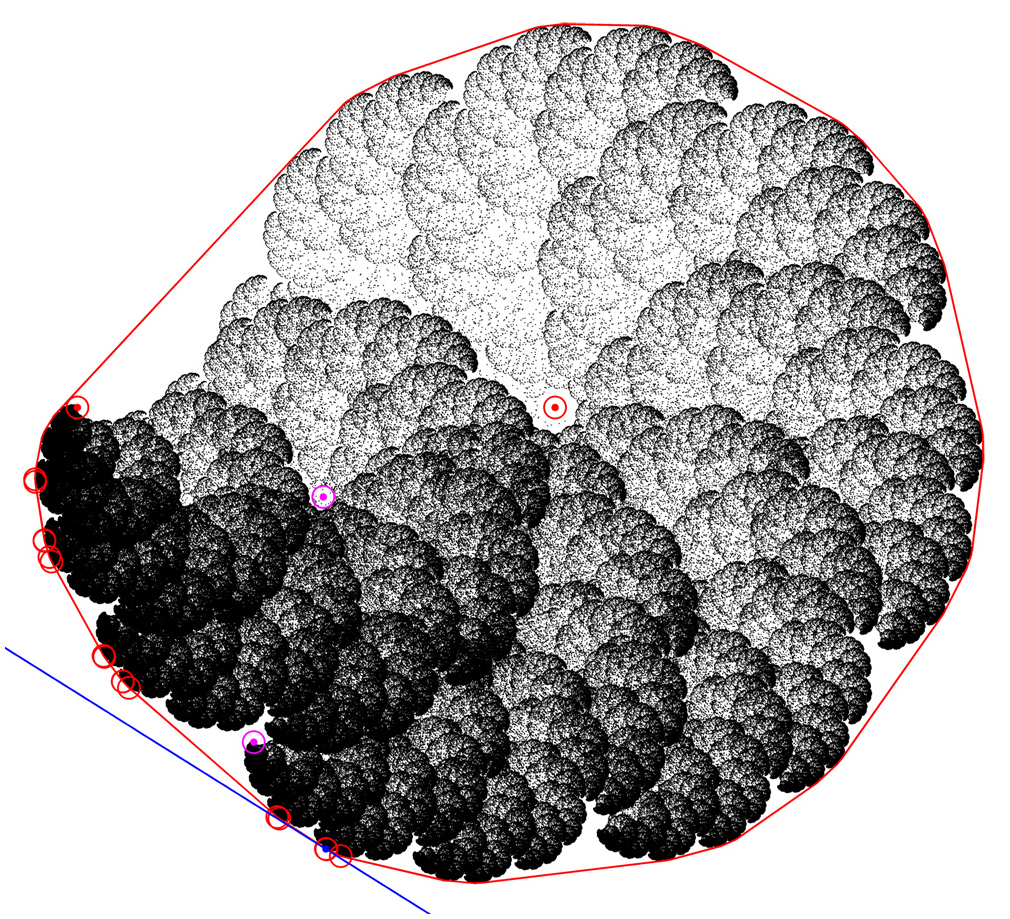}
\caption{A random C-IFS fractal.}
\label{s050409}
\end{figure}
\begin{figure}[H]
\centering
\includegraphics[width=360pt]{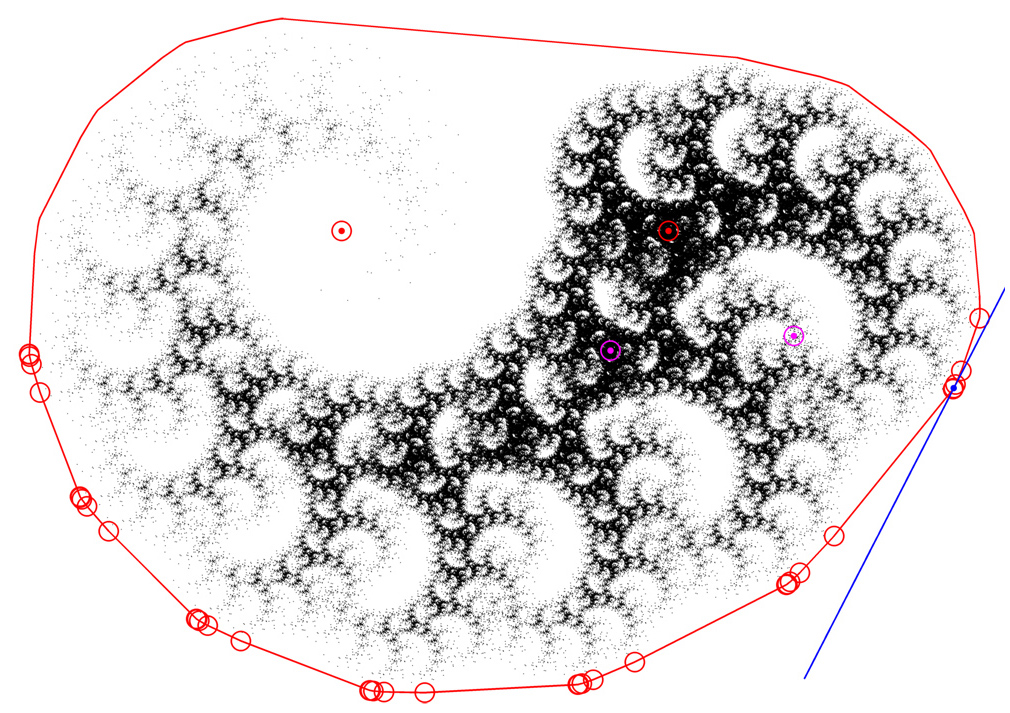}
\caption{Illustration of the method's predictive nature; plot iteration level $L=20$.}
\label{s050410}
\end{figure}
\begin{figure}[H]
\centering
\includegraphics[width=360pt]{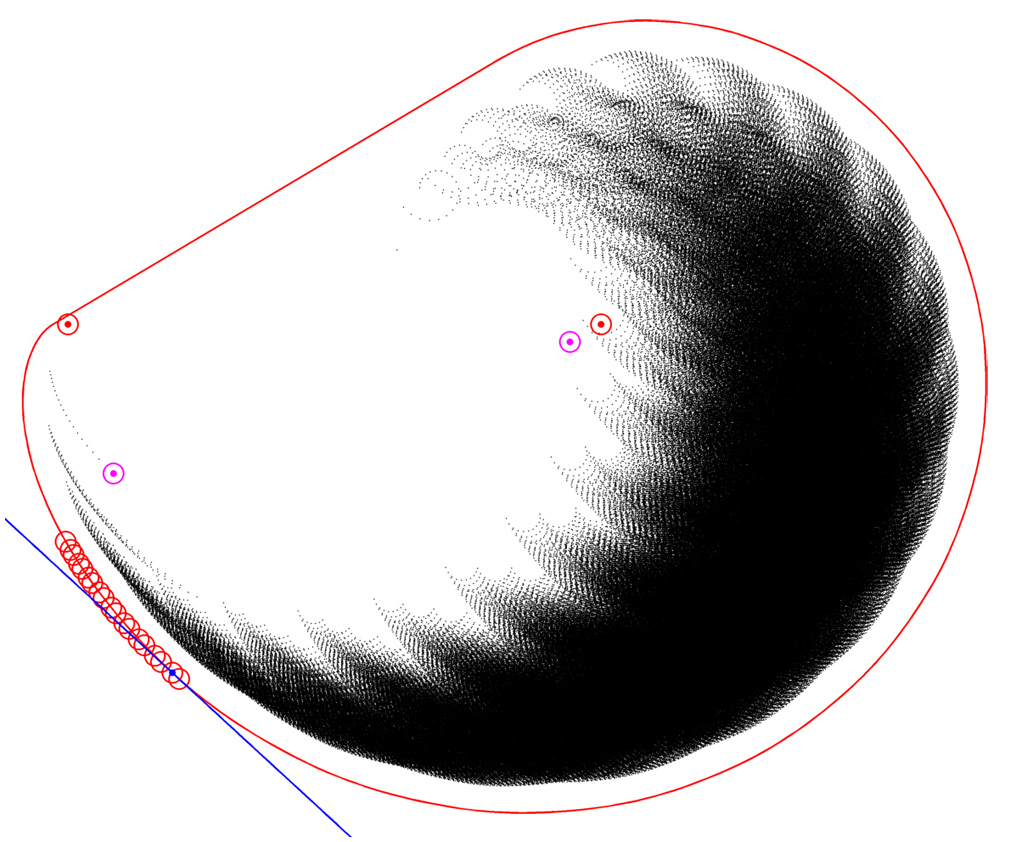}
\caption{Illustration of the method's predictive nature; plot iteration level $L=20$.}
\label{s050411}
\end{figure}

\section{Linear Optimization over IFS Fractals} \label{s0403}

An algorithm is described for maximizing a linear target function over a fractal of unity $F=\langle\mcT\rangle,\ \mcT=\{T_1,\ldots,T_n\}$ utilizing bounding circles and the Containment Lemma \ref{s020301}.

\begin{dfn} \label{s040301}
We say that a closed ball $C=B(c,r)=\{z\in\C: |z-c|\leq r\}$ with center $c\in\C$ and radius $r>0$ is a \textbf{bounding circle} to the fractal $F$ if it contains it $C\supset F$. It is an \textbf{ideal bounding circle} if it is invariant with respect to the Hutchinson operator $\mH_\mcT(C)\subset C$ and $c\in\mint\ \Conv(F)$.
\end{dfn}

The following $B(c,r)$ is an ideal bounding circle for any IFS fractal \cite{ba00135, ph00004}:
\[ c = \frac{1}{n}\sum_{k\in\mcN} p_k,\ \ r = \frac{\mu_*}{1-\vl_*}\max_{k\in\mcN} |p_k-c|\ \ \mathrm{where}\ \ \vl_*=\max_{k\in\mcN} |\cf_k|,\ \ \nu_*=\max_{k\in\mcN} |1-\cf_k|. \]
The method will search for maximizing truncation(s) recursively over the iterates $T_a(C)$. To make the search efficient, it compares and discards iterates which are ``dominated'' by others, with respect to a target vector $\tau$ and the corresponding target function.

\begin{dfn} \label{s040302}
We say that a finite address $a\in\mcAf$ \textbf{dominates} another $b\in\mcAf$ with respect to the target $\tau\in\C$, the IFS $\mcT=\{T_1,\ldots,T_n\}$, and the ideal bounding circle $C=B(c,r)\subset\C$, if $|a|=|b|$ and $\langle\tau,\ T_a(c)-T_b(c)\rangle \geq \vl_b r|\tau|$. Denoted as $a\succ_{\tau,\mcT,C} b$ or just $a\succ b$.

\vspace{0.3cm}
\begin{figure}[H]
\centering
\includegraphics[width=300pt]{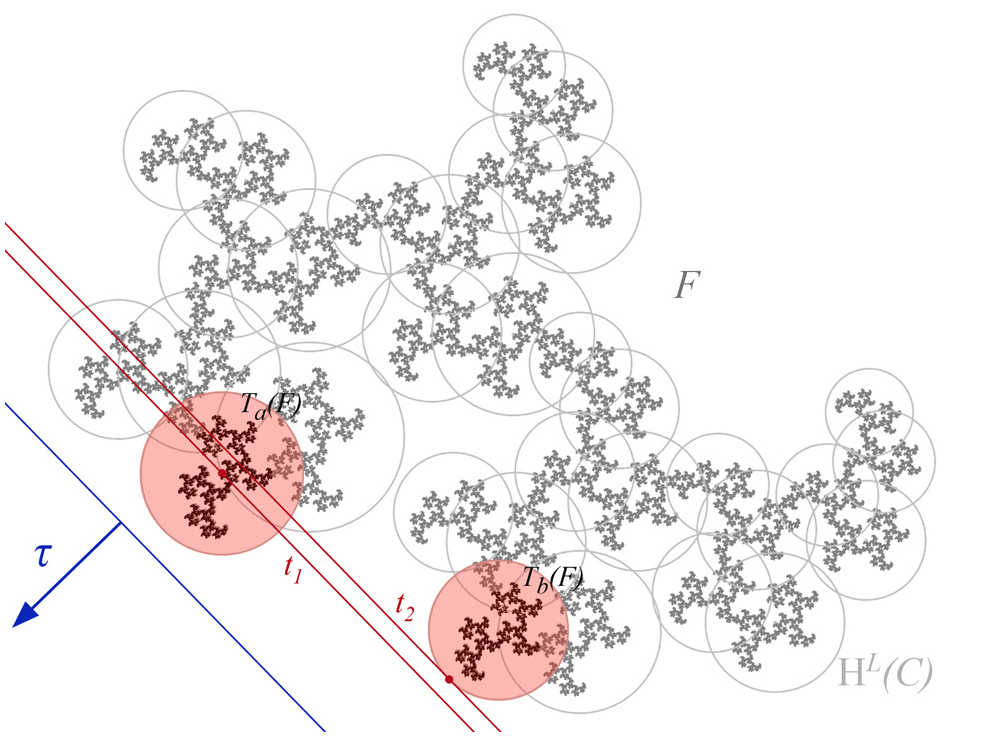}
\caption{Illustration of a finite address dominating another, for the parameters:\vspace{0.15cm}\ \newline \vspace{0.001cm} \hspace{4.35cm} $\cf_1=0.7\exp\left(\frac{5\pi}{15}\vi\right),\ \ \cf_2=0.6\exp\left(\frac{3\pi}{15}\vi\right)$.\vspace{0.15cm}\newline Note that the fixed point locations only affect scaling, not the geometry (see Section \ref{s0501}).}
\label{s040303}
\end{figure}
\end{dfn}

The inequality in the above definition is equivalent to the following
\[ t_1\libeq \langle\tau,\ T_a(c)\rangle \geq t_2\libeq \left\langle\tau,\ T_b(c)+\vl_b r\frac{\tau}{|\tau|}\right\rangle \]
which is illustrated on Figure \ref{s040303}. This implies that the centre of the iterated circle $T_a(C)\subset\mH^L(C), L\in\N$ has a greater-than-or-equal target value $t_1$ than any point (including the maximizing value $t_2$) over $T_b(C)\subset\mH^L(C)$. Due to $c\in\mint\ \Conv(F)$, we can infer that there must be a point in the subfractal $T_a(F)\subset \mH^L(F)=F$ which has a strictly larger target value than any of the points in the subfractal $T_b(F)\subset \mH^L(F)=F$, and therefore the maximizing algorithm can discard $T_b(F)$ and thus $T_b(C)$. This greatly increases the efficiency of seeking maximizers of $\tau$.

\begin{prp} \label{s040304}
The relation $\succ$ is a strict partial order over finite addresses, meaning it is an ordering relation that is irreflexive, transitive, and asymmetric.
\end{prp}
\noindent
\prf
We see that the relation is irreflexive, since $a\succ a$ would imply that $0=\langle\tau,\ T_a(c)-T_a(c)\rangle \geq \vl_b r|\tau| >0$ which is a contradiction.\\
To show transitivity, assuming that $a_1\succ a_2\succ a_3,\ a_{1,2,3}\in\mcAf$ we need $a_1\succ a_3$.
\[ \langle\tau,\ T_{a_1}(c)-T_{a_3}(c)\rangle = \langle\tau,\ T_{a_1}(c)-T_{a_2}(c)\rangle + \langle\tau,\ T_{a_2}(c)-T_{a_3}(c)\rangle \geq \]
\[ \geq \vl_{a_2}r|\tau| + \vl_{a_3}r|\tau| \geq \vl_{a_3}r|\tau|\ \ \rar\ \ a_1\succ a_3. \]
Lastly for asymmetry, we need that if $a\succ b$ then $b\succ a$ cannot hold. Clearly $a\succ b$ implies that $\langle\tau,\ T_a(c)-T_b(c)\rangle$ is strictly positive, but if $b\succ a$ also held then $\langle\tau,\ T_b(c)-T_a(c)\rangle= -\langle\tau,\ T_a(c)-T_b(c)\rangle$ would also be strictly positive, which is a contradiction. \sqr

\begin{dfn} \label{s040305}
We say that an element is \textbf{maximal} among the finite addresses $\mcA_0\subset\mcAf$ if no other element dominates it. Furthermore, denote the subset of maximal addresses as\ \ $\Argmax(\mcA_0)=\Argmax_{\tau,\mcT,C}(\mcA_0)\libeq \{a\in\mcA_0:\ \nexists b\in\mcA_0:b\succ_{\tau,\mcT,C} a\}$.
\begin{figure}[H]
\centering
\vspace{0.5cm}
\includegraphics[width=250pt]{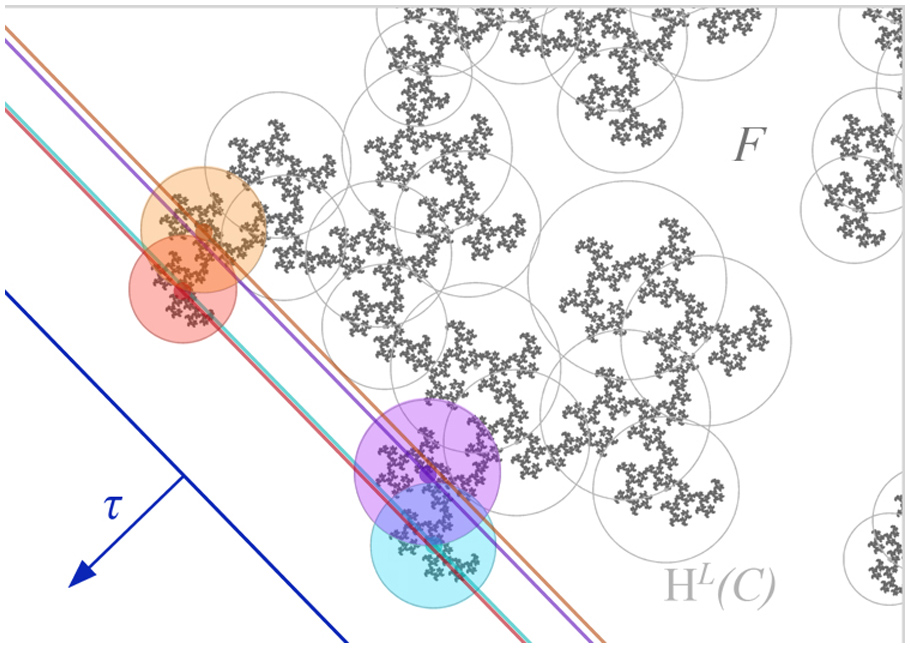}
\caption{Illustration of $\Argmax$ at some iteration level.}
\label{s040306}
\end{figure}
\end{dfn}

The above definition for the ``maximal element(s)'' in a subset of finite addresses is the standard way for partial ordering relations, illustrated in Figure \ref{s040306}. It depicts in the highlighted circle iterates the four maximal addresses with respect to the target $\tau$ and a bounding circle $C$, at some iteration level $L\in\N$ by the Hutchinson operator $\mH$. These addresses are considered to be maximal, since no other address dominates them at that iteration level. In the iterative maximization, we keep these addresses for the next iteration.

Finally, we arrive at the algorithm employing the above concepts. The algorithm uses the invariance of the ideal bounding circle $C\supset\mH(C)$ to eliminate the redundant non-maximal iterated circles of the form $T_a(C)$. The recursiveness lies in the step from potential maximizing truncations $\mcA_0$ to $\mcA_0\times\mcN$, corresponding to the subdivision of each circle $T_a(C),\ a\in\mcA_0$ into the sub-circles $T_a(T_k(C)),\ k\in\mcN$ to be compared for maximality. Only the maximal ones survive, and can again be subdivided iteratively for comparison. At each iteration
\[ \mamax_{f\in F}\ \langle\tau,f\rangle \subset \bigcup_{a\in\mcA_0} T_a(F) \subset \bigcup_{a\in\mcA_0} T_a(C). \]
So upon iteration, a set $\mcA_0\subset\mcAf$ of potential maximizing truncations are accumulated, and then the circles $T_a(C),\ a\in\mcA_0$ are subdivided further. According to Lemma \ref{s030203} and Theorem \ref{s030211}, as the subdivisions progress and the addresses in $\mcA_0$ get longer, all addresses in $\mcA_0\subset\mcAf$ eventually become blowable (see Definition \ref{s040201}), since extremal irreducible forms cannot be longer than $2|\nu(\mcAf)|$ according to Theorem \ref{s030211}. The blow-ups can then be verified to maximize the target function $z\mapsto\langle\tau,z\rangle\ (z\in\C)$. So the stopping criterion of the algorithm will be whether $\exists\blo(\mcA_0)$ holds, which must eventually occur within $2|\nu(\mcAf)|$ iterations. The returned output will be
\[ \ifrm\left(\mamax_{f\in\blo(\mcA_0)}\ \langle\tau,f\rangle\right) \]
representing the irreducible forms belonging to the maximizing blow-ups of $\mcA_0$. A target is maximized in either one or two extrema in $\blo(\mcA_0)\cap\Ext(F)$, so the above $\mamax$ set has either one or two elements.

\begin{algorithm}
\caption{(Linear Optimization Algorithm for IFS Fractals of Unity; LOAF)\newline Determines the irreducible forms of the maximizing extrema of a target direction $\tau\in\C$ over a fractal of unity $F=\langle\mcT\rangle$, using an ideal bounding circle $C=B(c,r)$. Call with $\mcA_0=\mcN$.} \label{s040307}
\begin{algorithmic}[1]
\Function{FractalLinOpt}{$\mcA_0,\tau,\mcT,C$}
    \If{$\exists\blo(\mcA_0)$}
        \State \textbf{return} $\ifrm\left(\mamax_{f\in\blo(\mcA_0)}\ \langle\tau,f\rangle\right)$
    \Else
        \State $\mcA_0\libeq \Argmax_{\tau,\mcT,C}(\mcA_0\times\mcN)$\Comment{Subdivision of $\mH^L(C)$ into $\mH^L(\mH(C))$.}
        \State \textbf{return} \Call{FractalLinOpt}{$\mcA_0,\tau,\mcT,C$}
    \EndIf
\EndFunction
\end{algorithmic}
\end{algorithm}

To ensure that the pseudocode can be implemented in a more efficient way as a program, we relate the elements of $\mcA_0$ to the calculation of $\Argmax_{\tau,\mcT,C}(\mcA_0\times\mcN)$. This corresponds to the subdivision of circles into their local iterates, and finding the maximal among them. Thus the calculation of $\Argmax$ is simplified via the following equivalences:
\[ a\in\Argmax_{\tau,\mcT,C}(\mcA_0\times\mcN)\ \lrar\ \nexists b\in\mcA_0\times\mcN: b\succ a\ \lrar \]
\[ \lrar\ \forall b\in\mcA_0\times\mcN:\ \left\langle\tau,\ T_a(c)+\vl_a r\frac{\tau}{|\tau|}\right\rangle > \langle\tau,\ T_b(c)\rangle\ \lrar \]
\[ \left\langle\tau,\ T_a(c)+\vl_a r\frac{\tau}{|\tau|}\right\rangle > \max_{b\in\mcA_0\times\mcN} \langle\tau,\ T_b(c)\rangle\ \ \mathrm{for\ each}\ \ a\in\mcA_0\times\mcN. \]
So we see that the next-level iterates of $C$ can be tested for maximality amongst one another, by first calculating the constant $\max$ in the last inequality, and then comparing it to each new target value on the left, for $a\in\mcA_0\times\mcN$. In order to make the above calculations even more efficient, we can keep track of the iterated centers $T_a(c)$ and the iterated fixed points $T_a(p_j),\ j\in\mcN$, since they imply the next-level iterates as follows:
\[ T_aT_k(c) = T_a(p_k)+\cf_k(T_a(c)-T_a(p_k))\ \ (a\in\mcA_0,\ k\in\mcN) \]
\[ T_aT_k(p_j) = T_a(p_k)+\cf_k(T_a(p_j)-T_a(p_k))\ \ (a\in\mcA_0,\ k,j\in\mcN). \]
These equations follow from the fact that $T_a$ is an affine map (Corollary \ref{s020304}).

\end{appendices}

\bibliographystyle{unsrt}
\bibliography{mybib}
\addcontentsline{toc}{section}{\textbf{References}}

\end{document}